\documentclass[12pt]{article}

\usepackage{amsmath,epsfig,epsf,psfrag,natbib}
\usepackage{amssymb}
\usepackage[latin1]{inputenc}
\usepackage{rotating}
\usepackage{amsbsy, url}
\usepackage{lscape, epstopdf}
\usepackage{color,latexsym,amsfonts,amsthm,amscd,graphicx}

\newcommand{\ben}{\begin{enumerate}}
\newcommand{\een}{\end{enumerate}}
\newcommand{\be}{\begin{equation}}
\newcommand{\ee}{\end{equation}}
\newcommand{\bas}{\begin{eqnarray*}}
\newcommand{\eas}{\end{eqnarray*}}
\newcommand{\ba}{\begin{eqnarray}}
\newcommand{\ea}{\end{eqnarray}}

\newtheorem{theorem}{Theorem}
\newtheorem{lemma}{Lemma}

\newtheorem{remark}{Remark}

\newtheorem{assumption}{Assumption}

\newtheorem{prop}{Proposition}

\newcommand{\e}{ { \mathbb{E}}}
\newcommand{\var}{ {\mathbb{V}\rm ar }}
\newcommand{\cov}{ {\mathbb{C}\rm ov}}

\newcommand{\pr}{{\rm pr}}

\def\T{{ \mathrm{\scriptscriptstyle \top} }}
\newcommand{\convergeto}{ {\overset{d}{\longrightarrow \; }}}

\usepackage{authblk}
\usepackage{algorithm}
\usepackage{algorithmic}

\begin{document}
\date{}
\title{Tuning-parameter-free optimal propensity
score matching approach for causal inference}

\author[1]{Yukun Liu\thanks{Corresponding author:  ykliu@sfs.ecnu.edu.cn}}
\author[2]{Jing Qin}
\affil[1]{
    KLATASDS-MOE,
 School of Statistics,
  East China Normal University,
    Shanghai 200062, China}
\affil[2]{ National Institute of Allergy and Infectious Diseases,
National Institutes of Health, Bethesda,  Maryland 20892, U.S.A.}
\renewcommand*{\Affilfont}{\small }
\renewcommand\Authands{ and }
\date{}
\maketitle

\begin{abstract}
Propensity score matching (PSM) is a pseudo-experimental method
that uses statistical techniques to construct an artificial
control group by matching each treated unit
with one or more untreated units of similar characteristics.
To date, the problem of determining the optimal number of matches per unit,
which plays an important role in PSM, has not been adequately addressed.
We propose a tuning-parameter-free PSM method
based on the nonparametric maximum-likelihood estimation
of the propensity score under the monotonicity constraint.
The estimated propensity score is piecewise constant,
and therefore automatically groups data. Hence,
our proposal is free of tuning parameters.
The proposed estimator is asymptotically
semiparametric efficient for the univariate case,
and achieves this level of efficiency in the multivariate case when the outcome and the propensity score depend on the covariate
in the same direction.
We conclude that matching methods based on the propensity score alone cannot, in general, be efficient.

\vspace{1cm}

\end{abstract}

\section{Introduction}

To assess the treatment effect in medical studies,
randomized and controlled clinical
trials are the gold standard because
the baseline covariates are balanced in the treatment and
control arms by the randomization. To evaluate the effectiveness of
an economic program or policy in econometrics or political
science, however, randomization is difficult or
impossible to implement for various reasons. In observational studies,
the available covariate information from people
who participated in the program or not may be unbalanced.
The simple two-sample $t$-test is likely to produce biased results.
It is desirable to replicate a randomized experiment as
closely as possible by obtaining treated and control
groups with similar covariate distributions. Due to
the simplicity and intuitiveness of adjusting the distribution of covariates among
samples from different populations, matching methods
are widely used in applied statistics,
econometrics, and epidemiology. Many examples
can be found in a comprehensive review paper by \cite{Stuart2010}.

A common feature of matching methods is the outcome-independence of matching, i.e., outcome values are not used in the
matching process even if they are available at the time
of matching.
This precludes the selection of a matched sample
that leads to a desired or undesired result.
Commonly used matching methods are all based on covariate
values. They use
the ``distance'' between treated and untreated individuals,
which is a measure of
the similarity between two individuals.
A well-known example is the Mahalanobis distance.
Another popular matching method is propensity matching,
which was proposed by \cite{Rosenbaum1983}.
The propensity score is the conditional probability of assignment to a
treatment given a vector of covariates.
Suppose that adjusting for a set of covariates is
sufficient to eliminate confounding.
A key observation made by \cite{Rosenbaum1983} is that adjusting for the propensity
score is also sufficient to eliminate confounding.
In contrast to covariate-based distance matching
methods, propensity score matching (PSM) has the
advantage of reducing the dimensionality of
matching to a single dimension,
making the matching process much easier.

Although various matching methods have
been used, the theoretical results
have only recently been studied by \cite{Abadie2006,Abadie2012,Abadie2016}. This series of papers showed that matching estimators based on
the covariate distance involve a biased term that
is only negligible under certain
regularity conditions.
Moreover, the matching estimators
are not necessarily root-$n$-consistent,
and some strong conditions are needed
to guarantee the root-$n$-consistency.
Furthermore, they demonstrated that, even in settings where
matching estimators are root-$n$-consistent,
simple matching estimators with a fixed number of
matches do not attain the semiparametric efficiency bound.
Furthermore, they found that matching estimators based on
the estimated propensity score have smaller variances than
those based on the true PSM.
To make the matching methods accessible to practitioners,
\cite{Imbens2015} used three examples to demonstrate practical implementations from the theoretical literature, and
provided detailed recommendations on how the procedures should be performed.
In her review paper in Statistical Science,
\cite{Stuart2010} lists some of the major software
packages that implement matching procedures. A regularly updated version
is also available at
\url{http://www.biostat.jhsph.edu/~stuart/propensityscoresoftware.html}.

In the PSM approach,
the most commonly used strategy is to divide the
interval $[0,1]$ into $K$ subintervals
$[a_i,a_{i+1}], i=1,2,...,K$, where $a_1=0$ and $ a_{K+1}=1$.
Then, one may compare the treatments and controls
for propensity scores that fall into the same subinterval.
The performance of this approach
is always influenced by the
choice of $K$, which is generally artificial,
and the numbers of individuals falling into each subinterval,
which is hard to control. Too many subintervals may
lead to inflated variances, whereas too few subintervals
may lead to biased results.
Existing matching methods
\citep{Abadie2006,Abadie2012,Abadie2016} require the
matched number of controls for each treatment case to be
fixed (not divergent to infinity).
A natural question is ``Are there any optimal methods
for choosing $K$ and the minimum number of
observations in each matched interval?''
To the best of our knowledge,
no research has yet been conducted
to address these issues.

In this paper, we present a tuning-parameter-free optimal matching method
to solve the two thorny issues outlined above,
where ``optimal'' means the resulting matching estimator
achieves the asymptotic semiparametric efficiency lower bound.
Rather than the commonly used parametric logistic model,
we assume a nonparametric and
monotone nondecreasing function for the propensity score
in the univariate case.
In the multivariate case,
we assume a single-index propensity score model
with a nonparametric and monotone nondecreasing link function.
Our method is based on the semiparametric maximum-likelihood
estimation of the propensity score function, which
can be implemented numerically by the well-known
pool adjacent violated algorithm \citep[PAVA]{Ayer1955}.
As the semiparametric maximum-likelihood estimator (MLE)
of the propensity score is a piecewise step function,
individuals in the treatment arm can be exactly
matched by individuals in the control arm based
on their estimated propensity scores.
This matching method is purely data-driven and
involves no artificial interference.
We show that the proposed tuning-free matching estimator
is not only asymptotically unbiased, but also achieves
the semiparametric efficiency lower bound in the univariate case.
For the multivariate case,
the proposed estimator also achieves the semiparametric efficiency
lower bound if the outcome and propensity score
depend on the covariate in the same direction.
Otherwise, the proposed estimator remains consistent, but is not efficient.
In this situation, other PSM
methods are relatively inefficient.
Our theoretical results depend critically on
the shape-restricted inference and empirical process theory.
Our numerical simulation results show that
the proposed method outperforms
existing commonly used PSM methods in terms of the mean square error.
A real econometric dataset is analyzed for illustration.
All technical proofs are given in the supplementary material for clarity.

\section{Efficient estimation under shape constraints}

\subsection{Setup}
We adopt the potential outcome framework
\citep{Neyman1923,Rubin1974}
with a binary treatment.
Let $Y(1)$ and $Y(0)$ be the potential outcomes of a treatment and a control,
respectively, which cannot be observed simultaneously.
Let $X$ be the baseline covariate,
and $D$ be the treatment indicator with $D=1$ denoting a treatment and $D$=0 denoting a control.
The outcome is $Y(1)$ if $D=1$ and $Y(0)$ otherwise, which
can be written as $Y =Y(D) = DY(1)+(1-D)Y(0)$.
Let $(Y_i, X_i, D_i), i=1,2,....,n$, be $n$ independent and
identically distributed observations from $(Y, X, D)$.
We focus on the estimation of
the average treatment effect on the treated (ATT)
\[
\tau = \e\{ Y(1) - Y(0)|D=1\};
\]
the average treatment effect $\e\{ Y(1) - Y(0)\}$ can be estimated similarly.
For the identifiability of treatment effects, we make the commonly used
unconfoundedness assumption
\citep{Rubin1978,Rosenbaum1983},
i.e., conditional on the observed covariates,
the treatment indicator is independent of the potential outcomes.

\begin{assumption}[Unconfounded Treatment Assignment]
\label{unconfoundness}
\(
D\perp (Y(0), Y(1)) |X.
\)
\end{assumption}

Denote the propensity score as $\pi(x)= \pr(D=1|X=x)$.
Under Assumption \ref{unconfoundness},
\cite{Rosenbaum1983} showed that $D$ and $(Y(0),Y(1))$
are conditionally independent
when $X$ is replaced by $\pi(X)$ in the condition, namely
\[
D\perp (Y(0), Y(1)) |\pi(X).
\]
The most commonly used model for the propensity score is
the logistic regression model;
however, misspecification of parametric models may lead to
inconsistent or misleading treatment effect estimators.
To alleviate the risk, we relax
the parametric model assumption to a
completely nonparametric model in the univariate case,
or a semiparametric index model in the multivariate case.

\begin{assumption}
\label{covariate-continuity}
The covariate $X$ has an absolutely continuous density and $\var(X)$ is positive
in the univariate case and positive-definite in the multivariate case.
In addition, the functions $\mu_0(X)= \e\{ Y(0)|X \}$,
$\mu_1(X)= \e\{ Y(1)|X \}$, $\sigma_0^2(X)=\var(Y(0)|X)$,
and $\sigma_1^2(X)=\var(Y(1)|X)$ are all well-defined for $P_X$-almost surely,
and the quantities $\e\{ Y(1)^2 \}, \e\{ Y(0)^2 \}$, $\e\{ \mu_1^2(X)\}$, and $\e\{ \mu_0^2(X)\}$
are all finite.
\end{assumption}

Assumption \ref{covariate-continuity}, which requires the covariate to be nondegenerate and
the outcome and covariate variables to have finite variances, is trivial.
Below, we consider the estimation problem of $\mu_1$ and $\tau$
in the univariate and multivariate cases separately.

\subsection{Estimation in the univariate case }

We do not make a parametric assumption on the propensity score, but instead
assume that $\pi(x)$ is a nonparametric and monotone nondecreasing
function.
Without loss of generality, we assume
\(
X_1\leq X_2\leq\ldots\leq X_n
\).
Based on $\{(X_i, D_i): 1\leq i\leq n\}$,
the likelihood of $\pi$ is a binomial likelihood of the form
\[
L_B (\pi) =\prod_{i=1}^n\pi(X_i)^{D_i}\{1-\pi(X_i)\}^{1-D_i}, \quad {\rm s.t.}\quad
\pi(X_1)\leq \pi(X_2)\ldots\leq \pi(X_n).
\]
By Theorem 2.12 of \cite{Barlow1972},
maximizing this likelihood with respect to $\pi$
is equivalent to minimizing
\[
\sum_{i=1}^n\{D_i-\pi(X_i)\}^2
\]
under the same monotonicity constraint.
The solution or the MLE $\hat \pi$ is a step function and
is the left derivative of the greatest convex minorant of the cumulative sum diagram
\cite[Theorem 1.1]{Barlow1972}.

Write $\hat{\pi}_i=\hat{\pi}(X_i)$, and let
$0=i_0<i_1<\cdots<i_k=n$ be the locations of the inflection points of
the greatest convex minorant of the cumulative sum diagram.
Then,
\ba
\label{ties}
\hat{\pi}_i=\hat{\pi}(X_i)=\hat{\pi}_{i_j}, \quad i_{j-1}<i\leq i_{j}, \quad j=1, \ldots, k.
\ea
According to the lemma on page 34 of \cite{Barlow1972},
\ba
\label{pi-ij}
\hat{\pi}_{i_j}=\frac{\sum_{l=i_{j-1}+1}^{i_{j}}D_l}{i_{j }-i_{j-1}}, \quad 1\leq j\leq k.
\ea
We propose to estimate $\mu_1 = \e\{Y(1)\}$ by
\ba
\label{proposed-mu1}
\hat{\mu}_1
= \frac{1}{n}\sum_{i=1}^n\frac{D_iY_i}{\hat{\pi}(X_i )}
= \frac{1}{n}\sum_{j=1}^k\frac{1}{\hat{\pi}_{i_j}}\sum_{l=i_{j-1}+1}^{i_{j}}D_lY_l
= \sum_{j=1}^k\rho_{j}\hat{\mu}_{1i_j},
\ea
where
\(
\rho_{j}=(i_{j}-i_{j-1})/n
\)
is the proportion of observations falling in the $j$-th interval $(x_{i_{j-1}},x_{i_{j}}]$, and
\(
\hat{\mu}_{1j}= \sum_{l=i_{j-1}+1}^{i_{j}}D_lY_l/\sum_{l=i_{j-1}+1}^{i_{j}}D_l
\)
is the group mean.
If $\hat{\pi}(X_i ) = 0$, the accompanying $D_i$ must also be zero,
and we define $D_i/ \hat{\pi}(X_i ) = 0/0$ as 0.
Essentially, $\hat{\mu}_1$ is a weighted average of subgroup means,
where the subgroups are formed by the steps
of the shape-restricted nonparametric MLE $\hat{\pi}$.
Note that this grouping method is
{\it automatically data-driven} and is {\it free from any tuning parameter.}

\begin{assumption}
\label{overlap-univariate}
(i) There exists $c_0>0$ such that $ c_0 \leq \pi(t) \leq 1-c_0$ for
all $t \in \mathcal{X}$.
(ii) $\pi(x)$ has a continuous derivative $\pi'(x)$ and
there exists $c_1>0$ such that $1/c_1\leq \pi'(x) \leq c_1$.
\end{assumption}

\begin{assumption}
\label{Lipschitz-mu1-univariate}
There exists $L>0$ such that
$|\mu_1(x_1) - \mu_1(x_2)|\leq L|x_1-x_2|$.
\end{assumption}

Assumption \ref{overlap-univariate} (i) is the commonly used overlap assumption
in the literature of causal inference. Under Assumption \ref{overlap-univariate} (ii),
both $\pi(x)$ and its inverse are Lipschitz-continuous.
Assumption \ref{overlap-univariate} (ii) and Assumption \ref{Lipschitz-mu1-univariate} imply that
the function $\mu_1(\pi^{-1}(t))$ is Lipschitz-continuous, which plays an important role
in governing random fluctuations in our proof.
Under these assumptions,
we show that the proposed estimator $\hat \mu_1$ is asymptotically unbiased, normal,
and achieves the semiparametric efficiency lower bound.

\begin{theorem}
\label{thm-mu1}
Suppose that
Assumptions \ref{unconfoundness}--\ref{Lipschitz-mu1-univariate}
are all satisfied.
The proposed estimator
$\hat{\mu}_1$ is asymptotically normal, i.e.,
$
\sqrt{n}(\hat{\mu}_1-\mu_1)\convergeto N(0,\sigma_\mu^2)
$ as $n\rightarrow\infty$,
where $\convergeto $ denotes ``converges in distribution to'' and
$
\sigma_\mu^2
=
\e\{ \sigma_1^2(X) / \pi(X) \}
+
\var\{ \mu_1(X) \}
$,
and the asymptotic variance attains the semiparametric efficiency lower bound.
\end{theorem}

By PSM, we propose to estimate $\tau$,
the ATT, by
\[
\hat{\tau} = \frac{1}{n_1} \sum_{i=1}^n D_i\left \{Y_i-\frac{\sum_{j=1}^n(1-D_j)Y_jI(\hat{\pi}(X_j)=\hat{\pi}(X_i))}
{\sum_{j=1}^n(1-D_j)I(\hat{\pi}(X_j)=\hat{\pi}(X_i))}
\right\},
\]
where $n_1=\sum_{i=1}^nD_i$.
A more concise form of $\hat{\tau} $ is provided in Lemma \ref{att-lemma}.

\begin{lemma}
\label{att-lemma}
The proposed PSM estimator $\hat{\tau}$ can be equivalently expressed as
\bas
\hat{\tau} &=&
\frac{1}{n_1}\sum_{j=1}^n\left \{ D_jY_j-(1-D_j)Y_j\frac{\hat{\pi}_j}{1-\hat{\pi}_j}\right \},
\eas
which is an inverse probability weighting estimator.
\end{lemma}

With the estimated propensity scores $\hat \pi_i$,
the estimator of $\tau $ developed by \cite{Hirano2003} is
\bas
\tilde{\tau} =\frac{1}
{\sum_{j=1}^n\hat{\pi}_j} \sum_{i=1}^n \left\{ D_iY_i-(1-D_i)Y_i \frac{ \hat{\pi}_i}{1-\hat{\pi}_i} \right\}.
\eas
By Theorem 1.7 of \cite{Barlow1972},
the shape-restricted MLEs $\hat{\pi}_i $ satisfy $\sum_{i=1}^n(\hat{\pi}_i-D_i)=0$
or, equivalently,
$\sum_{i=1}^n\hat{\pi}_i=\sum_{i=1}^nD_i=n_1$.
We find that $\tilde \tau = \tilde \tau$, i.e.,
the proposed PSM estimator is equal to
that of \cite{Hirano2003} in the form for
the ATT.

\begin{assumption}
\label{Lipschitz-mu0-univariate}
There exists $L>0$ such that
$|\mu_0(x_1) - \mu_0(x_2)|\leq L|x_1-x_2|$.
\end{assumption}

Assumption \ref{Lipschitz-mu0-univariate} and
Assumption \ref{overlap-univariate} (ii) imply that
the function $\mu_0(\pi^{-1}(t))$ is Lipschitz-continuous, which
is used to govern random fluctuations in our proof of
the asymptotic normality of $\hat \tau$.

\begin{theorem}
\label{thm-tau}
Suppose that
Assumptions \ref{unconfoundness}--\ref{overlap-univariate} and \ref{Lipschitz-mu0-univariate}
are all satisfied.
Let $\eta = \pr(D=1)$ and $\tau(X)= \mu_1(X)- \mu_0(X)$.
The proposed PSM estimator for the ATT is asymptotically normal, i.e.,
$
\sqrt{n}(\hat{\tau} -\tau ) \convergeto N(0,\sigma_{\tau}^2)
$ as $n\rightarrow\infty$,
where
\[
\sigma_{\tau}^2=
\frac{1}{\eta^2}\e\left [ \pi (X)\{ \tau(X)-\tau \}^2+\pi(X)\sigma_1^2(X)+\frac{\pi^2(X)}{1-\pi(X)}\sigma_0^2(X)\right ],
\]
and the asymptotic variance attains the semiparametric efficiency lower bound.
\end{theorem}

\cite{Hahn1998} derived the semiparametric efficiency lower bound
for the estimation of $\tau$ when the propensity score is unknown (his Theorem 1)
and when it is known (his Theorem 2).
Our asymptotic variance $\sigma_{\tau}^2$ is exactly equal to
the asymptotic semiparametric efficiency lower bound when the propensity score is unknown.
Therefore, it attains the semiparametric efficiency lower bound.

\section{Estimation in the multivariate covariate case}

When the covariate is a $d$-variate ($d>1$),
we assume that the propensity score is
\ba
\label{multi-para-ps}
\pr(D=1|X=x)=\pi(x^\T\beta),
\ea
where $\pi$ is a monotone nondecreasing function
and $\beta$ is an unknown true $d$-variate parameter.
We assume that $ \|\beta\| = 1$ for identifiability.
Suppose that
a consistent estimator $\hat \beta$ of $\beta$ is available.
We shall consider two estimation methods for $\beta$ in subsections 3.1 and 3.2.
Let
\(
Z_i(\hat \beta)=X_i^\T\hat{\beta}
\)
for $i=1,2,\ldots,n$.
Without loss of generality, we assume that
\(
Z_1(\hat \beta)\leq Z_2(\hat \beta)\leq \cdots\leq Z_n(\hat \beta)
\).
The log-likelihood of $\pi$ becomes
\[
\sum_{i=1}^n \left[ D_i\log\pi(Z_i(\hat \beta))+(1-D_i)\log\{1-\pi(Z_i(\hat \beta))\}\right].
\]
The MLE of $\pi$ is a step function determined by
$\pi(Z_1(\hat \beta)), \pi(Z_2(\hat \beta)), \ldots \pi(Z_n(\hat \beta))$,
which satisfy the monotonicity restriction
\(
\pi(Z_1(\hat \beta))\leq \pi(Z_2(\hat \beta))\cdots\leq \pi(Z_n(\hat \beta)).
\)
Denote the shape-restricted MLE of $\pi$
as $\hat{\pi}(\cdot)$.
Our proposed PSM estimators for $\mu_1$ and $\tau$ are
\ba
\hat{\mu}_1
&=&
\frac{1}{n}\sum_{i=1}^n\frac{D_iY_i}{\hat{\pi}(Z_i(\hat \beta))}
=\frac{1}{n}\sum_{i=1}^n\frac{D_iY_i}{\hat{\pi}(X_i\hat{\beta})},
\label{est-mu-multi-general}
\\
\hat{\tau}
&=&
\frac{1}{n_1}\sum_{j=1}^n\left \{ D_jY_j-(1-D_j)Y_j\frac{\hat{\pi}(X_j\hat{\beta})}{1-\hat{\pi}(X_j\hat{\beta})}\right \},
\label{est-tau-multi-general}
\ea
respectively, where we have used
a multivariate version of Lemma \ref{att-lemma}.

\begin{assumption}
\label{range-XB}
The ranges of $X$ and $\beta$,
$\mathcal{X}$ and $\mathcal{B}$, are convex and compact.
Let
$t_{\rm low} = \inf\{ X^\T \gamma: X \in \mathcal{X}, \gamma \in \mathcal{B}\}-\varepsilon_0$
and
$t_{\rm up} = \sup\{ X^\T \gamma: X \in \mathcal{X}, \gamma \in \mathcal{B}\}+\varepsilon_0$
for some $\varepsilon_0>0$.

\end{assumption}

\begin{assumption}
\label{pi-multi}
There exists $c_0\in (0, 1)$ such that
$ c_0 \leq \pi(t) \leq 1-c_0$
and $\pi$ has a continuous second derivative on
$ [t_{\rm low}, \ t_{\rm up}]$,
where
$t_{\rm low} $
and
$t_{\rm up}$ are defined in Assumption \ref{range-XB}.

\end{assumption}

Under Assumption \ref{pi-multi},
$\pi'(t)$ is also continuous on the
closed interval $ [t_{\rm low}, \ t_{\rm up}]$.
Therefore, it must be Lipschitz-continuous, i.e.,
there exists $c_1>0$ such that
$ | \pi'(t) - \pi'(s)| \leq c_1 |s-t|$ for any $t_{\rm low}\leq s, t\leq t_{\rm up}$.

\begin{assumption}
\label{Assumption-X-Gamma}
There exists a constant $M>0$ such that
the density function $f_{X^\T \gamma}(u)$ of
$X^\T \gamma$ satisfies
$ f_{X^\T \gamma}(u) \leq M$ for all $x\in \mathcal{X}$
and $\gamma \in \mathcal{B}$.

\end{assumption}

Define
$
\mu_{1}^* (u; \gamma)=\e\{ Y(1)|X^\T \gamma = u \} = \e\{ \mu_1(X)|X^\T \gamma = u\},
$
and
$
\mu_{0}^* (u; \gamma)=\e\{ Y(0)|X^\T \gamma = u \} = \e\{ \mu_0(X)|X^\T \gamma = u\}
$.

\begin{assumption}
\label{Lipschitz-mu1-multivariate}
The function $\mu_{1}^* (u; \gamma)$ is continuous in both $u$ and $\gamma$.
\end{assumption}

\begin{assumption}
\label{Lipschitz-mu0-multivariate}
The function $\mu_{0}^* (u; \gamma)$ is continuous in both $u$ and $\gamma$.
\end{assumption}

Because $\mathcal{X}$ and $\mathcal{B}$ are both compact,
Assumptions \ref{Assumption-X-Gamma} and \ref{Lipschitz-mu1-multivariate}
imply that the function $\mu_{1}^* (X^\T \gamma_1; \gamma_2)$ is Lipschitz-continuous
with respect to $(\gamma_1, \gamma_2)$, i.e.,
there exists a constant $L$ such that
\bas
|
\mu_{1}^* (X^\T \gamma_1; \gamma_2) - \mu_{1}^* (X^\T \gamma_3; \gamma_4)
|
\leq
L (\|\gamma_1-\gamma_3\| + \|\gamma_2 - \gamma_4\|), \;
\gamma_1, \gamma_2, \gamma_3, \gamma_4 \in \mathcal{B}.
\eas
The function $\mu_0^*$ has the same property under
Assumptions \ref{Assumption-X-Gamma} and \ref{Lipschitz-mu0-multivariate}.
In general, if $\hat\beta$ is $\sqrt{n}$-consistent and asymptotically normal,
the proposed estimators for $\mu_1$ and $\tau$
both follow asymptotically normal distributions,
and both are asymptotically semiparametric efficient
under certain additional conditions.
Let $\mathbb{P}_n$ denote the empirical measure based on data $\{(X_i, D_i, Y_i): 1\leq i\leq n\}$.

\begin{theorem}
\label{mu-multi-general-ps}
Suppose that model \eqref{multi-para-ps} is true,
and that
Assumptions \ref{unconfoundness}--\ref{covariate-continuity}
and \ref{range-XB}-\ref{Lipschitz-mu1-multivariate}
are satisfied.
Define
\bas
B_1 = \e \left[
\frac{\pi'(X^\T \beta) }{\pi(X^\T {\beta})}
\{ \mu_1(X) - \mu_1^* ( X^\T \beta; \beta) \} X^\T
\right].
\eas
If $\hat \beta - \beta = O_p(n^{-1/2})$, then the following results hold as $n\rightarrow\infty$.
\ben
\item[(1)]
A linear approximation for $\hat \mu_1$ is
\ba
\hat \mu_1 -\mu_1
&=&
\mathbb{P}_n \Big[
\frac{D- \pi(X^\T \beta)}{\pi(X^\T \beta)}
\{ \mu_1(X ) - \mu_1^* ( X^\T \beta; \beta ) \}
\nonumber \\
&&
+
\frac{D (Y -\mu_1(X ))}{\pi (X^\T \beta )}+\mu_1(X ) -\mu_1 \Big]
+
B_1 (\hat \beta - \beta) + o_p(n^{-1/2}).
\label{mu-approx-general}
\ea

\item[(2)]
If $\mu_1(X) = \tilde \mu_1(X^\T \beta)$ for some function $\tilde \mu_1(\cdot)$, then
$\mu_1^* ( X^\T \beta; \beta) = \tilde \mu_1 ( X^\T \beta )$ and
\bas
\sqrt{n}(\hat \mu_1 -\mu_1)
&=&
\sqrt{n} \mathbb{P}_n \left \{\frac{D (Y -\tilde \mu_1(X^\T \beta))}{\pi (X^\T \beta )}+\tilde \mu_1(X^\T \beta) -\mu_1\right \}
+ o_p(1) \\
&
\convergeto& N(0, \sigma_{\mu, \rm m}^2),
\eas
where
$
\sigma_{\mu, \rm m}^2
= \var(\tilde \mu_1(X^\T \beta) ) + \e \{\sigma_1(X)( 1- \pi (X^\T \beta ))/\pi (X^\T \beta ) \}
$
is the asymptotic semiparametric efficiency lower bound.
Namely, if both the propensity score and regression function
depend on the covariate $X$ in the same direction,
then $\hat \mu_1$ achieves the asymptotic semiparametric
efficiency lower bound, which holds for any $n^{1/3}$-consistent estimator $\hat \beta$.

\een

\end{theorem}

\begin{theorem}
\label{tau-multi-general-ps}
Suppose that model \eqref{multi-para-ps} is true,
and that
Assumptions \ref{unconfoundness}--\ref{covariate-continuity},
\ref{range-XB}--\ref{Assumption-X-Gamma}, and \ref{Lipschitz-mu0-multivariate}
are satisfied. Define
\ba
\label{B-2}
B_2
=
\e \left[ \frac{ \pi' (X^\T \beta) }{1-\pi(X^\T \beta)}
\{ \mu_0(X ) - \mu_0^* ( X^\T {\beta}; \beta) \} X^\T \right].
\ea
If $\hat \beta - \beta = O_p(n^{-1/2})$, then the following results hold as $n\rightarrow\infty$.
\ben
\item[(1)]
A linear approximation for $\hat \tau$ is
\ba
\hat{\tau} - \tau
&=&
\frac{1}{\eta}
\mathbb{P}_n\{ D (\tau(X) - \tau) + D (Y(1) - \mu_1(X))\}
\nonumber \\
&&
-
\frac{1}{\eta}
\mathbb{P}_n \left[ (1-D )\{ Y(0) - \mu_0(X)\} \frac{ {\pi} (X^\T \beta)}{1- {\pi}(X^\T \beta)}
\right] \nonumber \\
&&
+
\frac{1}{\eta} \mathbb{P}_n \left[ \frac{ D- \pi (X^\T \beta) }{1- \pi (X^\T \beta)}
\{ \mu_0(X ) - \mu_0^* ( X^\T {\beta}; \beta) \} \right]\nonumber
\\
&&
-
\frac{1}{\eta} B_2 (\hat \beta - \beta) + o_p(n^{-1/2}).
\label{tau-approx-general}
\ea

\item[(2)]
If $\mu_0(X) = \tilde \mu_0(X^\T \beta)$ for some function $\tilde \mu_0(\cdot)$, then
$\mu_0(X) = \mu_0^* ( X^\T \beta; \beta) = \tilde \mu_0 ( X^\T \beta )$ and
\bas
\sqrt{n}(\hat{\tau} - \tau)
&=&
\frac{1}{\eta} \sqrt{n} \mathbb{P}_n
\Big[ D (\tau(X) - \tau) + D (Y(1) - \mu_1(X)) \\
&&
\hspace{2cm}
-
(1-D )\{ Y(0) - \mu_0(X)\} \frac{ {\pi} (X^\T \beta)}{1- {\pi}(X^\T \beta)} \Big]
+ o_p(1) \\
& \convergeto& N(0, \sigma_{\tau, \rm m}^2),
\eas
where
\bas
\sigma_{\tau, \rm m}^2
&=&
\frac{1}{\eta^2} \Big[\e\{ \pi(X) (\tau(X) - \tau )^2 \}
+ \e( \pi(X^\T \beta) \sigma_1^2(X) )
+
\e( \sigma_0^2(X) \frac{ \{\pi(X^\T \beta)\}^2}{ 1- \pi(X^\T\beta) } ) \Big]
\eas
is the asymptotic semiparametric efficiency lower bound (Theorem 1 of Hahn, 1998).
Namely if both the propensity score and regression function
depend on covariate $X$ in the same direction,
then $\hat \tau$ achieves the asymptotic semiparametric
efficiency lower bound, which holds for any $n^{1/3}$-consistent estimator $\hat \beta$.

\een

\end{theorem}

Besides the asymptotic normality and efficiency results,
Theorems \ref{mu-multi-general-ps} and \ref{tau-multi-general-ps}
also indicate that
if the propensity score and regression functions
depend on the covariate $X$ in different directions,
or the regression functions do not obey single-index models,
then neither $\hat \mu_1$ nor $\hat \tau$ is asymptotically semiparametric efficient.

\begin{remark}
\cite{Imai2014} introduced a covariate balancing propensity score methodology
that models treatment assignment while
optimizing the covariate balance.
Suppose that $\pi(X^\T\beta)$
is a correctly specified model for the propensity score.
Observing the fact that, for any function $h(X)$,
\[
\e\left \{ \frac{D}{\pi(X^\T\beta)}h(X)-\frac{1-D}{1-\pi(X^\T\beta)}h(X)\right \}=0,
\]
instead of estimating $\beta$ by the maximum-likelihood method, \cite{Imai2014}
proposed to estimate $\beta$ by solving
\[
\mathbb{P}_n \left[ \frac{D - \pi(X^\T \beta)}{\pi(X^\T\beta)\{1-\pi(X^\T\beta)\}}h(X) \right] =0
\]
for some function $h(x)$. For example, $h(X)=X$, $\pi'_{\beta}(X\beta) X$, or the vector
consisting of all the linear and quadratic terms of $X$.
They argue that if the propensity score model is misspecified,
the MLE of
the propensity score might not balance the covariates, while
their proposed approach can balance the first and second
moments between the two arms.
If the dimension of $h(X)$ is greater than that of $\beta$,
this is a well-known over-identified estimation problem.
They estimate $\beta$ using the generalized method of moments \citep{Hansen1982}
and the empirical likelihood method \citep{Qin1994}.
The idea of matching with empirical likelihood was also considered by \cite{Qin2007}
for handling missing data problems.
In general, a higher dimension of $h$
produces more efficient estimators,
but creates a heavier computational burden.
In practical applications, one has to make a trade-off between computational cost
and estimation efficiency.

Let $\hat \pi(\cdot)$ be the MLE
of $\pi(\cdot)$ under the monotonicity constraint
based on observations $\{X_i^\T \hat \beta: 1\leq i\leq n\}$
for any given $\hat \beta$.
By the characterization of such a shape-restricted MLE \citep{Barlow1972},
we have
\[
\mathbb{P}_n [ \{ D - \hat{\pi}(X^\T \hat \beta) \} h( \hat{\pi}(X^\T \hat \beta) ) ] =0
\]
for any function $h$. Moreover, we can show that
\[
\mathbb{P}_n [ \{ D -\hat{\pi}(X^\T \hat{\beta}) \} h(X^\T \hat{\beta}) ] =o_p(n^{-1/2})
\]
if $h$ and $\pi$ are sufficiently smooth (see the proof of Lemma 9 in the supplementary material).
In other words, our proposal can balance any covariate function of the form $h(X\hat{\beta})$
up to a higher than root-$n$ asymptotic order $o_p(n^{-1/2})$.

\end{remark}

\subsection{When $\pi(\cdot)$ is known}

We first consider a simple case in which the function $\pi$
is known. For such cases, we could pretend that $\pi$ is unknown and monotone increasing,
and apply the proposed estimation procedure.
It is then natural to estimate $\beta$ by its MLE
$
\hat \beta = \arg\max \ell_B(\beta),
$
where
\[
\ell_B(\beta)=\prod_{i=1}^n [ D_i \log\{ \pi(X_i^\T\beta)\} + (1-D_i) \log \{ 1-\pi(X_i^\T\beta) \}].
\]

\begin{assumption}
\label{known-pi}
\ben
\item[(1)]
The true parameter value
$\beta $ is an interior point of $\mathcal{B}$
and
the unique solution to
\bas
\e
\left[
\frac{D -\pi(X^\T \tilde \beta)}{ \pi(X^\T \tilde \beta) \{ 1- \pi(X^\T \tilde \beta)\} }
\pi'(X^\T \tilde \beta) X \right] = 0
\eas
with respect to $\tilde \beta$.

\item[(2)]
The matrix $B_3$ is nonsingular, where
\bas
B_3 =
\e
\left[
\frac{ \{ \pi'(X^\T \beta ) \}^2}{ \pi (X^\T \beta ) \{ 1-\pi (X^\T \beta ) \}}
X X^\T \right].
\eas
\een
\end{assumption}

If Assumptions \ref{range-XB}, \ref{pi-multi}, and \ref{known-pi} are all satisfied,
then $\hat \beta$ is $\sqrt{n}$-consistent to $\beta$ and admits a linear approximation
(see Lemma 11 in the supplementary material).

\begin{theorem}
\label{thm-multi-para-ps}
Suppose that model \eqref{multi-para-ps} is true,
$\pi$ is a known and monotone nondecreasing function,
and that
Assumptions \ref{unconfoundness}--\ref{covariate-continuity} and
\ref{range-XB}--\ref{known-pi}
are satisfied.
Then, the following results hold as $n\rightarrow\infty$.
\ben
\item[(1)]
A linear approximation for $\hat \mu_1$ is
\bas
\sqrt{n}(\hat \mu_1 -\mu_1)
&=&
\sqrt{n}\mathbb{P}_n \left[
\frac{D- \pi(X^\T \beta)}{\pi(X^\T \beta)}
\left\{ \mu_1(X) - \mu_1^* ( X^\T \beta; \beta) +\frac{ B_1 B_3^{-1} \pi'(X^\T \beta) X }{ 1- \pi(X^\T \beta) }
\right\}
\right. \\
&&
+
\left. \frac{D (Y(1) -\mu_1(X ))}{\pi (X^\T \beta )}+\mu_1(X) -\mu_1\right]
+ o_p(n^{-1/2}) \\
&
\convergeto& N(0, \sigma_{\mu, \rm kn}^2),
\eas
where
\bas
\sigma_{\mu, \rm kn}^2
&=& \var( \mu_1(X) ) + \e \left \{\frac{ 1- \pi (X^\T \beta ) }{ \pi (X^\T \beta ) } \sigma_1(X) \right \} \\
&&
+
\e \left[
\frac{1- \pi(X^\T \beta)}{\pi(X^\T \beta)}
\left\{ \mu_1(X) - \mu_1^* ( X^\T \beta; \beta) +\frac{ B_1 B_3^{-1} \pi'(X^\T \beta) X }{ 1- \pi(X^\T \beta) }
\right\}^2
\right].
\eas
\item[(2)]
A linear approximation for $\hat \tau$ is
\bas
\sqrt{n}(\hat{\tau} - \tau)
&=&
\sqrt{n} \frac{1}{\eta} \mathbb{P}_n
\Big[ D (\tau(X) - \tau) + D (Y(1) - \mu_1(X)) \\
&&
-
(1-D )\{ Y(0) - \mu_0(X)\} \frac{ {\pi} (X^\T \beta)}{1- {\pi}(X^\T \beta)} \\
&&
+
\frac{ D- \pi (X^\T \beta) }{1- \pi (X^\T \beta)}
\Big\{ \mu_0(X ) - \mu_0^* ( X^\T {\beta}; \beta)
- \frac{ \pi'(X^\T \beta) }{ \pi(X^\T \beta)} B_2 B_3^{-1} X
\Big\} \Big] \nonumber \\
&& + o_p(n^{-1/2}) \\
&
\convergeto&
N(0, \sigma_{\tau,\rm kn}^2),
\eas
where
\bas
\sigma_{\tau, \rm kn}^2
&=&
\frac{1}{\eta^2} \e\Big[ \pi(X^\T \beta) (\tau(X) - \tau )^2 + \pi(X^\T \beta) \sigma_1^2(X)
+
\sigma_0^2(X) \frac{ \{\pi(X^\T \beta)\}^2}{ 1- \pi(X^\T\beta) } \\
&&
+
\frac{ \pi (X^\T \beta) }{1- \pi (X^\T \beta)}
\Big\{ \mu_0(X ) - \mu_0^* ( X^\T {\beta}; \beta)
- \frac{ \pi'(X^\T \beta) }{ \pi(X^\T \beta)} B_2 B_3^{-1} X
\Big\}^2
\\
&&
+
2
\Big\{ \pi (X^\T \beta) (\mu_0(X ) - \mu_0^* ( X^\T {\beta}; \beta))
- \pi'(X^\T \beta) B_2 B_3^{-1} X
\Big\} \cdot (\tau(X) - \tau )
\Big].
\eas
\een

\end{theorem}

Because the link function
$\pi $ is known, we may replace $\hat \pi(\cdot)$
in \eqref{est-mu-multi-general} and \eqref{est-tau-multi-general} by $\pi(\cdot)$ and consider the following two estimators:
\ba
\label{estimator-pava-m1-par}
\breve{\mu}_1
= \frac{1}{n}\sum_{i=1}^n\frac{D_iY_i}{ \pi (X_i\hat{\beta})},
\quad
\breve{\tau}
=
\frac{1}{n_1}\sum_{j=1}^n\left \{ D_jY_j-(1-D_j)Y_j\frac{\pi (X_j\hat{\beta})}{1-\pi (X_j\hat{\beta})}\right \}.
\ea

\begin{prop}
\label{thm-multi-para-ps-par}
Under the conditions of Theorem \ref{thm-multi-para-ps}, as $n\rightarrow\infty$,
$
\sqrt{n}(\breve{\mu}_1 - \mu_1 ) \convergeto N(0, \bar\sigma_{\mu, \rm kn}^2),
$
where
\bas
\bar \sigma_{\mu, \rm kn}^2
&=& \var( \mu_1(X ) ) + \e \left \{\frac{ 1- \pi (X^\T \beta ) }{ \pi (X^\T \beta ) } \sigma_1(X) \right \} \\
&&
+
\e \left[\frac{1 - \pi(X^\T \beta ) }{ {\pi}(X^\T \beta)}
\left\{\mu_1(X ) + \frac{ B_2 B_3^{-1} \pi'(X^\T \beta) X }{ 1- \pi(X^\T \beta) }
\right\}^2
\right]
\eas
and
\ba
\label{B-4}
B_4 =
\e \left\{
\frac{ \mu_1(X )}{ \pi(X^\T \beta) } \pi'(X^\T \beta)X^\T \right\},
\ea
and
$\sqrt{n}(\breve{\tau} - \tau) \convergeto N(0, \bar\sigma_{\tau, \rm kn}^2)$,
where
\bas
\bar \sigma_{\tau, \rm kn}^2
&=&
\frac{1}{\eta^2} \e \Big[
\pi (X^\T \beta ) (\tau(X)- \tau)^2 + \pi(X^\T \beta) \sigma_1^2(X)
+
\frac{\pi^2 (X^\T \beta )}{1-\pi (X^\T \beta )} \sigma_0^2(X)
\\
&&
+
\frac{ \pi (X^\T \beta )}{1-\pi (X^\T \beta )}
\Big\{ \mu_0(X)
-
\frac{ \pi'(X^\T \beta)}{ \pi(X^\T \beta)} B_4 B_3^{-1} X \Big\}^2 \\
&&
+2
(\tau(X)- \tau) \cdot
\Big\{ \mu_0(X) \pi(X^\T \beta)
-
\pi'(X^\T \beta) B_4 B_3^{-1} X \Big\}
\Big].
\eas
\end{prop}

\cite{Hahn1998} stated that, even if the propensity score is completely known,
the asymptotic semiparametric efficiency lower bound is the same as
that in the case where the propensity score is completely unknown.
Here, we assume a weaker assumption, namely
that the propensity score satisfies a single-index model with a known link function.
Therefore, the asymptotic semiparametric efficiency lower bound must also be
the same as that in the case where the propensity score is completely unknown.
When $Y(1)$ depends on $X$ in the same direction $\beta$ as
$D$ does, the proposed estimator $\hat \mu_1$, which
does not use the true link function $\pi(\cdot)$ but its PAVA estimator,
achieves the asymptotic semiparametric efficiency lower bound
$
\sigma_{\mu, \rm m}^2
= \var(\tilde \mu_1(X^\T \beta) ) + \e \{\sigma_1(X)( 1- \pi (X^\T \beta ))/\pi (X^\T \beta ) \}.
$
In contrast, the estimator $\breve{\mu}_1 $
has an asymptotic variance of
$
\bar \sigma_{\mu, \rm kn}^2
$, which is clearly greater than $\sigma_{\mu, \rm m}^2$.
In other words,
the estimator $\breve{\mu}_1 $ using the true link function $\pi(\cdot)$
is not asymptotically semiparametric efficient.
The estimation of $\tau$ has a similar property.
When $Y(0)$ depends on $X$ in the same direction $\beta$ as
$D$ does, the proposed estimator $\hat \tau$, which
does not use the true link function $\pi(\cdot)$ but its PAVA estimator,
achieves the asymptotic semiparametric efficiency lower bound
$\sigma_{\tau, \rm m}^2$.
However, the estimator $\breve{\tau} $ using the true link function $\pi(\cdot)$
does not.
These observations coincide with the results of \cite{Hahn1998}.

\subsection{When $\pi(\cdot)$ is unknown}

The known link function assumption may appear to be too strong,
as it is rarely known in practice.
In this subsection, we assume that the link function
$\pi(\cdot)$ in \eqref{multi-para-ps} is unknown and
monotone nondecreasing. We propose to estimate $\beta$ via the
simple score estimator (SSE) $\hat \beta$ of
\cite{Balabdaoui2019}.
Once $\hat \beta$ has been obtained, we estimate $\pi(\cdot)$
by PAVA based on the observations $(D_i, X_i^\T \hat \beta)$ ($1\leq i\leq n$).

We briefly review the SSE of \cite{Balabdaoui2019}.
We again use $\beta$ to denote the true index.
To identify $\beta$, we assume that $\|\beta\| = 1$ and that its
first nonzero component is positive.
Given $\gamma$, let $Z_i(\gamma) = X_i^\T \gamma$ and assume that
$Z_1(\gamma)\leq Z_2(\gamma) \leq\ldots\leq Z_n(\gamma)$. Let $\hat \pi_{\gamma}$
denote the PAVA estimator of $\pi(\cdot)$ that minimizes
$
\sum_{i=1}^n \{ D_i - \pi(Z_i(\gamma)) \}^2.
$
Define a $d-1$-dimensional sphere as $\mathcal{S}_{d-1} =
\{\gamma: \gamma \in \mathbb{R}^d, \|\gamma\| = 1\}$,
a one-to-one map $\mathbb{S}: [0, \pi]^{(d-2)}\times [0, 2\pi] \mapsto \mathcal{S}_{d-1}$ as
\ba
&&
\hspace{-0.4cm}
(\zeta_{(1)}, \zeta_{(2)}, \ldots, \zeta_{(d-1)})
\mapsto
(\cos(\zeta_{(1)}), \ \sin(\zeta_{(1)}) \cos(\zeta_{(2)}), \nonumber \\
&&
\hspace{3.9cm} \sin(\zeta_{(1)}) \sin(\zeta_{(2)}) \cos(\zeta_{(3)}), \ldots, \nonumber \\
&&
\hspace{3.9cm}
\sin(\zeta_{(1)}) \ldots \sin(\zeta_{(d-2)}) \cos(\zeta_{(d-1)}), \nonumber\\
&&
\hspace{3.9cm} \sin(\zeta_{(1)}) \ldots \sin(\zeta_{(d-2)}) \sin(\zeta_{(d-1)})),
\label{big-S-fun}
\ea
and a $d\times (d-1)$ matrix as
$J(\zeta) = \partial \mathbb{S}^\T (\zeta)/ \partial \zeta $.
Let $\zeta_0$ satisfy $\beta = \mathbb{S}(\zeta_0)$ and
$\hat \zeta$ be a {\it zero-crossing} of the function
\ba
\label{function-zero-crossing}
\phi_n(\zeta) = \mathbb{P}_n [ J^\T(\zeta) X\{D - \hat \pi_{\mathbb{S}(\zeta) }( X^\T \mathbb{S}(\zeta))\} ]
\ea
(see page 521 of \cite{Balabdaoui2019}
for the definition of {\it zero-crossing}).
Accordingly, we estimate $\beta$ by $\hat \beta = \mathbb{S}(\hat \zeta)$,
and estimate the propensity score function by $\hat \pi_{\hat \beta}(\cdot \hat \beta)$.
The resulting PSM estimators for $\mu_1$ and $\tau$ are
\bas
\hat{\mu}_1
=
\frac{1}{n}\sum_{i=1}^n\frac{D_iY_i}{\hat{\pi}_{\hat \beta}(X_i^\T\hat{\beta})}
\quad
\mbox{and}\quad
\hat{\tau}
=
\frac{1}{n_1}\sum_{j=1}^n\left \{ D_jY_j-(1-D_j)Y_j\frac{\hat{\pi}_{\hat \beta}(X_j^\T\hat{\beta})}{1-\hat{\pi}_{\hat \beta}(X_j^\T\hat{\beta})}\right \},
\eas
respectively.
To study the large-sample properties of the two estimators,
we assume the following conditions,
which correspond to Assumptions A3, A5, A7, and A9, respectively, of \cite{Balabdaoui2019}.

\begin{assumption}
\label{Ass3}
There exists $\delta_0>0$ such that the function
$
\pi_{ \gamma}(u) = \e\{ \pi(X^\T \beta)| X^\T \gamma =u\} $
is monotone increasing
on $I_{\gamma} = \{ X^\T \gamma: X\in \mathcal{X} \}$
for all $\gamma \in \mathcal{B}(\beta, \delta_0) = \{\gamma: \|\gamma - \beta\|\leq \delta_0\}$.

\end{assumption}

\begin{assumption}
\label{Ass5}
The distribution of $X$ admits a density $g$ that is differentiable
on $\mathcal{X}$. In addition, there exist positive constants
$c_1, c_2, c_3, c_4>0$ such that
$c_1\leq g\leq c_2$ and $c_3\leq \partial g/\partial x_j\leq c_4$
on $\mathcal{X}$ for all $1\leq j\leq d$.
\end{assumption}

\begin{assumption}
\label{Ass7}
For all $\gamma \neq \beta$ such that $\mathbb{S}(\gamma) \in \mathcal{B}(\beta, \delta_0)$,
the random variable
\[
\cov[(\zeta_0 - \zeta)^\T J^\T (\zeta) X, \
\pi(X^\T \mathbb(\zeta_0)) | X^\T \mathbb(\zeta)
]
\]
is almost surely not equal to zero.
\end{assumption}

\begin{assumption}
\label{Ass9}
$J^\T(\zeta_0) \e\{ \pi'(X^\T \beta) \cov(X|X^\T \beta) \} J(\beta)$
is nonsingular.
\end{assumption}

If Assumptions \ref{range-XB}, \ref{known-pi}(2), and \ref{Ass3}--\ref{Ass9} are satisfied, then
Theorem 3 of \cite{Balabdaoui2019}
implies that the estimator $\hat \beta = \mathbb{S}(\hat \zeta)$ is consistent and asymptotically normal
(see Lemma 12 in the supplementary material).

\begin{theorem}
\label{multi-nonpara-ps}
Suppose that model \eqref{multi-para-ps} is true, $\pi$ is an unknown and monotone nondecreasing function,
and that
Assumptions \ref{unconfoundness}--\ref{covariate-continuity},
\ref{range-XB}--\ref{Lipschitz-mu0-multivariate}, and \ref{Ass3}--\ref{Ass9}
are satisfied.
Define
\ba
\label{B-5}
B_5 = J (\zeta_0)
\{
J^\T(\zeta_0) \e[ \pi'(X^\T \beta) \var(X|X^\T \beta) ]J(\zeta_0)
\}^{-1} J^\T(\zeta_0).
\ea
Then, the following results hold as $n\rightarrow\infty$.
\ben
\item[(1)]

$
\sqrt{n} (\hat \mu_1 -\mu_1)
\convergeto N(0, \sigma^2_{\mu, \rm un}),
$
where
\bas
\sigma^2_{\mu, \rm un}
&=&
\e \left[
\frac{ 1- \pi(X^\T \beta) }{\pi(X^\T \beta) } \Big\{
\mu_1(X ) - \mu_1^* ( X^\T \beta; \beta )
+ B_1 B_5( X - \e(X|X^\T \beta ))\pi(X^\T \beta) \Big\}^2
\right] \\
&&
+
\e \left \{\frac{ \sigma_1^2(X ))}{\pi (X^\T \beta )} \right \}
+
\e \{ \mu_1(X ) -\mu_1 \}^2.
\eas

\item[(2)]
$
\sqrt{n}(\hat \tau - \tau)
\convergeto N(0, \sigma_{\tau, \rm un}^2),
$
where
\bas
\sigma_{\tau, \rm un}^2
&=&
\frac{1}{\eta^2}\e
\Big[ \pi(X^\T\beta) \sigma_1^2(X) \Big]
+
\frac{1}{\eta^2}\e\Big[ \sigma^2_0(X)\frac{ \{ \pi (X^\T \beta)\}^2 }{1- {\pi}(X^\T \beta)}
\Big] \\
&&
+
\frac{1}{\eta^2}\e
\{ \pi (X^\T \beta) (\tau(X) - \tau)^2 \}
+
\frac{1}{\eta^2}\e\Big[ \pi (X^\T \beta) \{1- \pi (X^\T \beta)\}\\
&&
\Big\{ \frac{ \mu_0(X ) - \mu_0^* ( X^\T {\beta}; \beta) }{1- \pi (X^\T \beta)}
-B_2 B_5( X - \e(X|X^\T \beta )) \Big\}^2 \Big] \\
&&
+
\frac{2}{\eta^2}\e
\Big[ (\tau(X) - \tau) \pi(X^\T \beta) \{1- \pi (X^\T \beta)\}
\Big\{ \frac{ \mu_0(X ) - \mu_0^* ( X^\T {\beta}; \beta) }{1- \pi (X^\T \beta)} \\
&&
-B_2 B_5( X - \e(X|X^\T \beta )) \Big\} \Big].
\eas

\een

\end{theorem}

Theorem \ref{multi-nonpara-ps} indicates that,
in general, if the link function $\pi(\cdot)$ is unknown and estimated by PAVA,
the proposed estimators $\hat \mu_1$ and $\tau$
are still consistent and asymptotically normal.
However, by Theorems \ref{mu-multi-general-ps} and \ref{tau-multi-general-ps},
they are not asymptotically semiparametric efficient
if $Y(1)$ or $Y(0)$ does not depend on the covariate in the same direction
as the treatment indicator does.

\section{Simulations}

To evaluate the finite-sample performance of the proposed estimators,
we conduct simulations to compare the following estimation methods:
\begin{itemize}
\item PAVA-MLE: the proposed PSM method with $\beta$ estimated by the MLE under the logistic propensity score model
and $\pi$ estimated by PAVA;
\item PAVA-SSE: the proposed PSM method with $\beta$ estimated by SSE and $\pi$ estimated by PAVA;
\item PARA: the proposed PSM estimator \eqref{est-tau-multi-general}
with $\hat \pi$ replaced by the logistic function and
$\hat \beta$ being the MLE under the logistic propensity score model;
\item PSM-$M$: the PSM method
with the propensity score estimated by the logistic regression model
and each case matched with $M$ controls.
Four choices of $M$ are considered: 3, 5, 10, and 15.
\end{itemize}

To generate data, we consider the bivariate $X=(X_1, X_2)$,
a linear logistic propensity score model
\bas
\pr(D=1|x_1,x_2)=\pi(2+x_1+x_2),
\eas
and the following regression models:
\bas
Y(1)
=
- ( X_1+X_2 )^a+\epsilon, \quad
Y(0)
=
3 h(X_1, X_2) -( X_1+ b X_2 )^a+\epsilon,
\eas
where $\epsilon, X_1$, and $ X_2 $ are independent and identically distributed as
$N(0, 1)$.
We choose $\pi(t) = e^t/(1+e^t)$ or the standard normal distribution function,
$ h(X_1, X_2) = \cos( X_1+ b X_2)$ (Model 1), $X_1$ (Model 2), $a=1, 2$, and $b=1, 0, -1$.
From each case, we generate 1000 samples with a sample size of $n =$ 500 and calculate
the seven estimators for the ATT $\tau$.

We first examine the results in Table \ref{logistic-tab}, which presents
the biases and root mean square errors (RMSEs) of the seven estimators
when $\pi(t) = e^t/(1+e^t)$.
The overall rate of nonmissing data is about $ \pr(D=1)= 81.6\% $.
As $\pi(t)$ is the logistic function, the seven estimators under comparison all have
correctly specified propensity score models.
The propensity score satisfies model \eqref{multi-para-ps} with
$\beta = (1/\sqrt{2}, 1/\sqrt{2})^\T$.
In all cases, although having negligible biases, the PARA estimator (which
uses the true logistic propensity score function)
always has the largest RMSE,
meaning that it is always the most unreliable
among the seven estimators under comparison.
This coincides with the finding of Hirano et al. (2003) that ``weighting by the
inverse of a nonparametric estimate of the propensity score, rather than the true propensity
score, leads to an efficient estimate of the average treatment effect.''

Under Model 1, the regression function in the control group is
a single-index model
$\mu_0(X) = \cos( \sqrt{1+b^2} \cdot X^\T \theta) - ( \sqrt{1+b^2} \cdot X^\T \theta )^a$
with $\theta = (1, b)/\sqrt{1+b^2}$.
When $b=1$, $\theta = \beta$. By Theorem \ref{tau-multi-general-ps},
the proposed estimator $\hat \tau$ in \eqref{est-tau-multi-general}
is asymptotically semiparametric efficient, regardless of
whether $\beta$ is estimated by the MLE or SSE.
We see from Table \ref{logistic-tab} that PAVA-SSE has
very similar performance to PAVA-MLE,
and both of them perform uniformly better than
the five competitors in terms of bias and RMSE.
The performance of the PSM estimator is dramatically influenced by
the number of matches, $M$, per unit; the PSM estimator has increasing
biases and RMSEs as $M$ increases from 3 to 15, and
PSM-15 with $M=15$ has twice the biases and RMSEs as PSM-3 with $M=3$.

When $b \neq 1$, we have $\theta \neq \beta$.
The proposed estimator $\hat \tau$ loses its semiparametric optimality.
Even so, when $b=0$ (the angle between $\theta$ and $\beta$ is 45 degrees),
the proposed estimators PAVA-SSE and PAVA-MLE
still achieve better performances than the other estimators,
but the relative advantage is smaller.
When $b= -1$, $\theta$ is perpendicular to $\beta$,
and the relative advantage of our estimators decreases further
until PSM becomes comparable or even better.
Surprisingly, we find that, in this case, the PSM estimator has decreasing
biases and RMSEs as $M$ increases from 3 to 15,
which is contrary to the case where $\theta = \beta$.
These findings indicate that the optimal choice of $M$ for the PSM method
critically depends on the true regression function $\mu_0(X)$,
and without any information on $\mu_0(X)$, it is impossible
to correctly specify the optimal $M$.
Additionally, the performance of the proposed estimator may be improved
by making use of information on $\mu_0(X)$.

Under model 2, $\mu_0(X)$ does not follow a single-index model,
and it cannot be written as $\tilde \mu_0(X^\T \theta)$
for some function $\tilde \mu_0$.
By Theorem \ref{tau-multi-general-ps},
the proposed estimator $\hat \tau$
is no longer semiparametric efficient,
but is still asymptotically normal.
The results in Table \ref{logistic-tab} corresponding to Model 2
suggest that, compared with the PSM estimators,
the proposed estimators are at least comparable
and perform uniformly better in some cases.

Table \ref{probit-tab} presents the results when $\pi(t) $ is chosen to be the standard normal distribution function.
In this situation, the overall rate of nonmissing data is about $ \pr(D=1)= 87.5\% $ and
only the proposed PAVA-SSE has a correctly specified propensity score model.
Compared with the four PSM estimators, the proposed PAVA-MLE and PAVA-SSE estimators
have uniformly smaller RMSEs and generally smaller biases.
Again, the two proposed estimators exhibit rather similar performance,
although
the index coefficient in the PAVA-MLE method
is estimated by the MLE under the logistic propensity score model.
As $M$ increases from 3 to 15, the PSM method may perform worse
in some cases and better in other cases.
When the true propensity score model changes from
a logistic model to a probit model,
the RMSEs of the PAVA-MLE and PAVA-SSE estimators
increase by no more than 26\% in eight out of the twelve cases.
In contrast, the RMSE of PSM-3 increases by at least 60\% in all cases,
and can be as large as 200\% (e.g., the case with Model 2, $a=1$, and $b=-1$).
This suggests that the PSM estimators are more sensitive
than the proposed estimator to the misspecification of
the propensity score model.

\begin{figure}[p]
\centering
\includegraphics[width=0.45\textwidth, height=0.24\textheight]{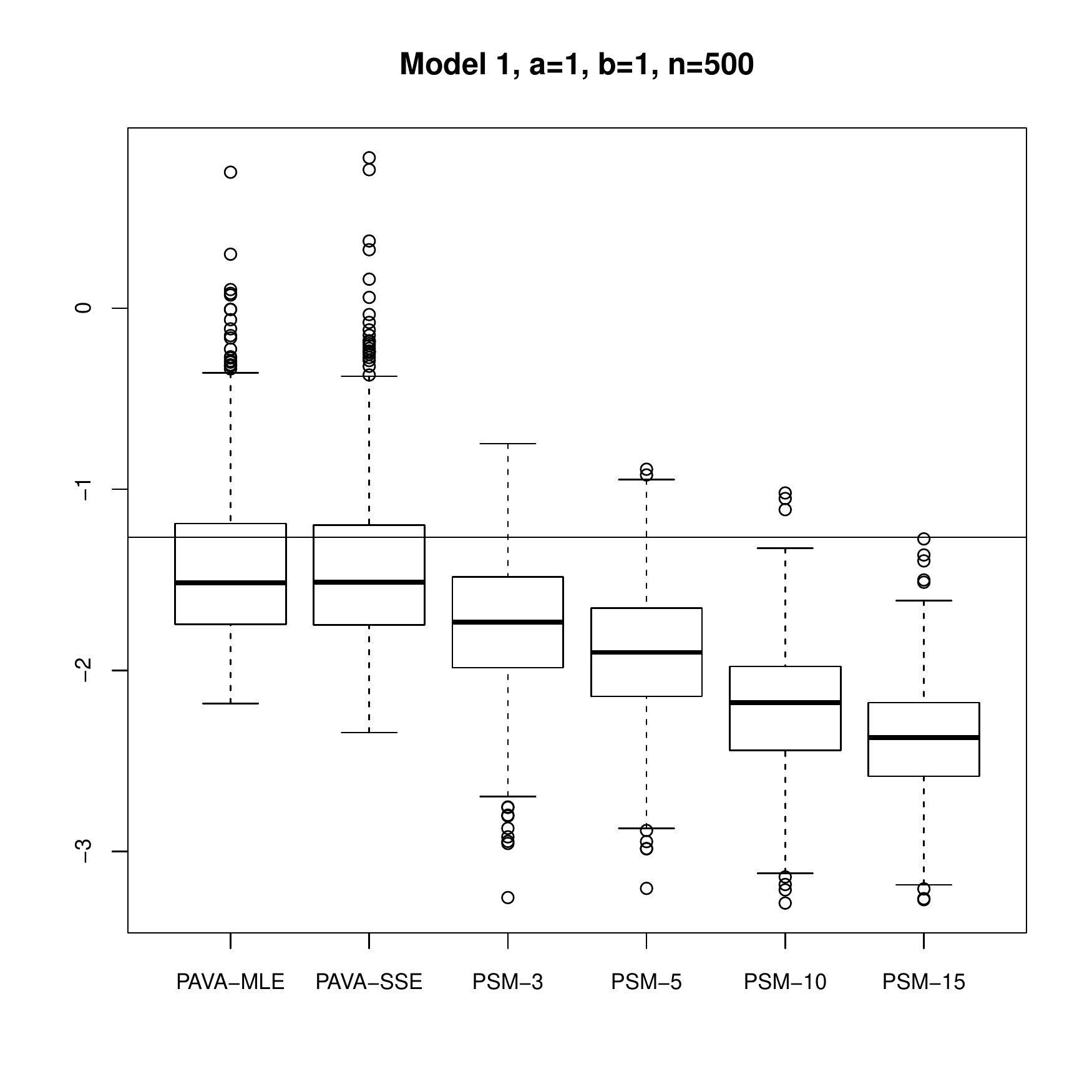}
\includegraphics[width=0.45\textwidth, height=0.24\textheight]{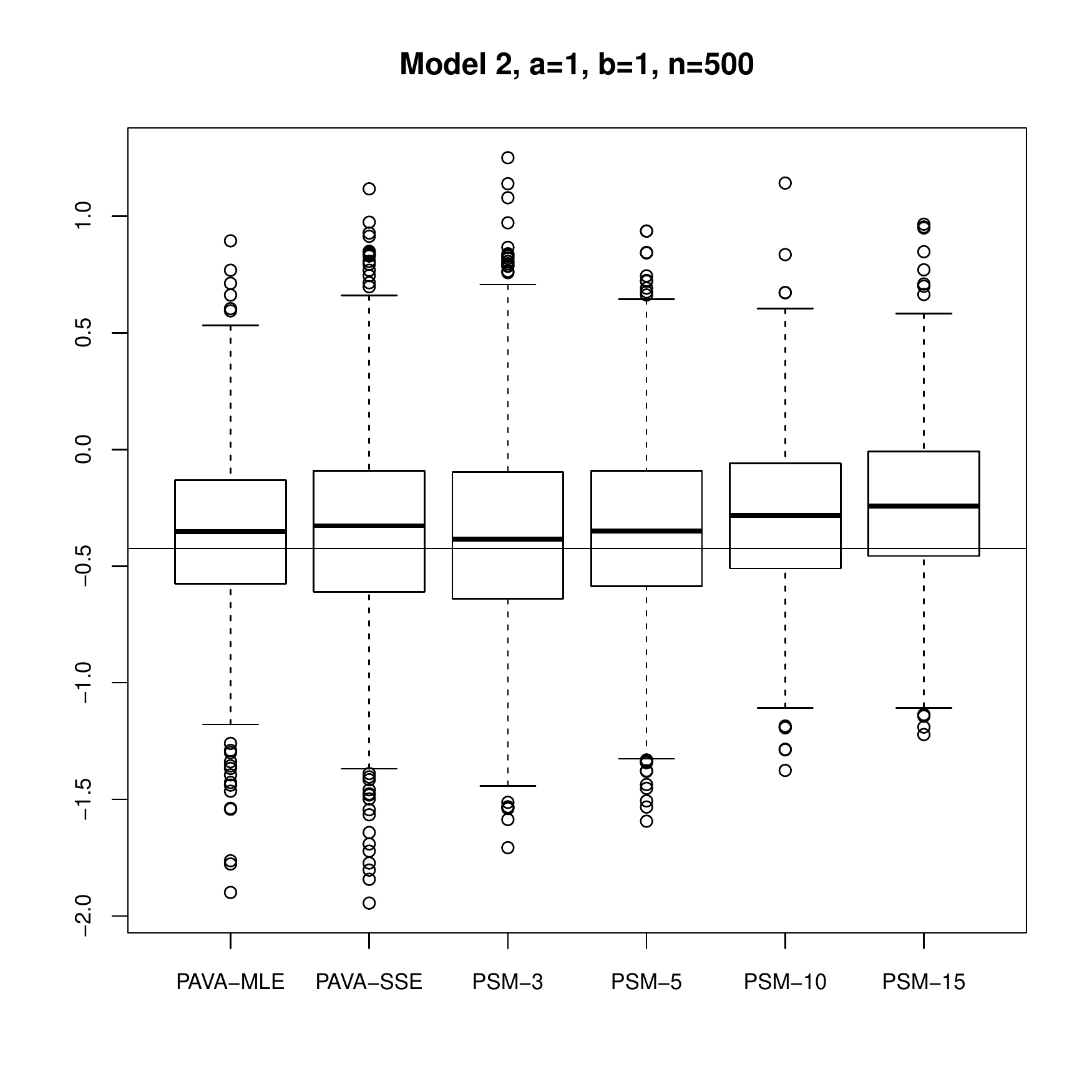} \\
\includegraphics[width=0.45\textwidth, height=0.24\textheight]{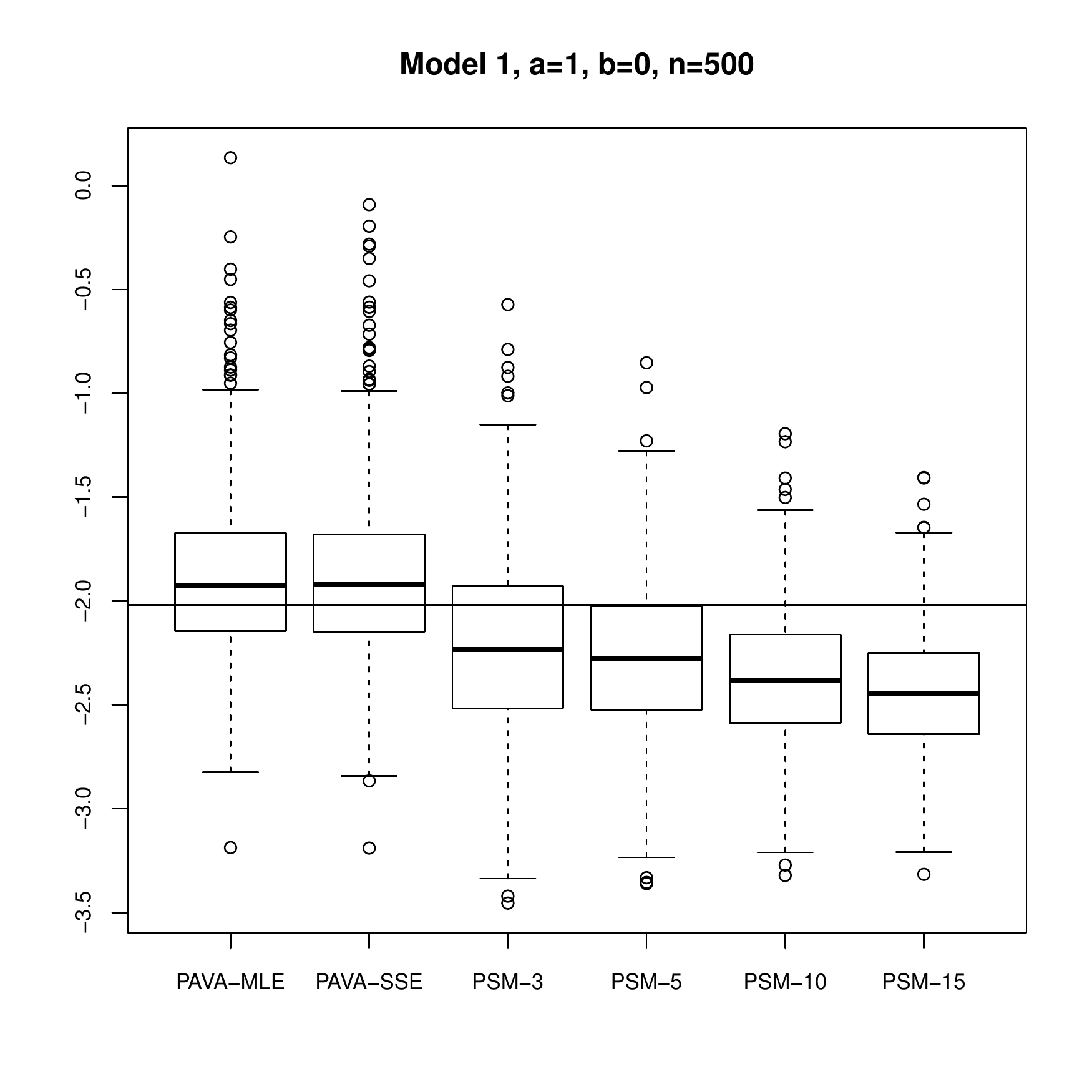}
\includegraphics[width=0.45\textwidth, height=0.24\textheight]{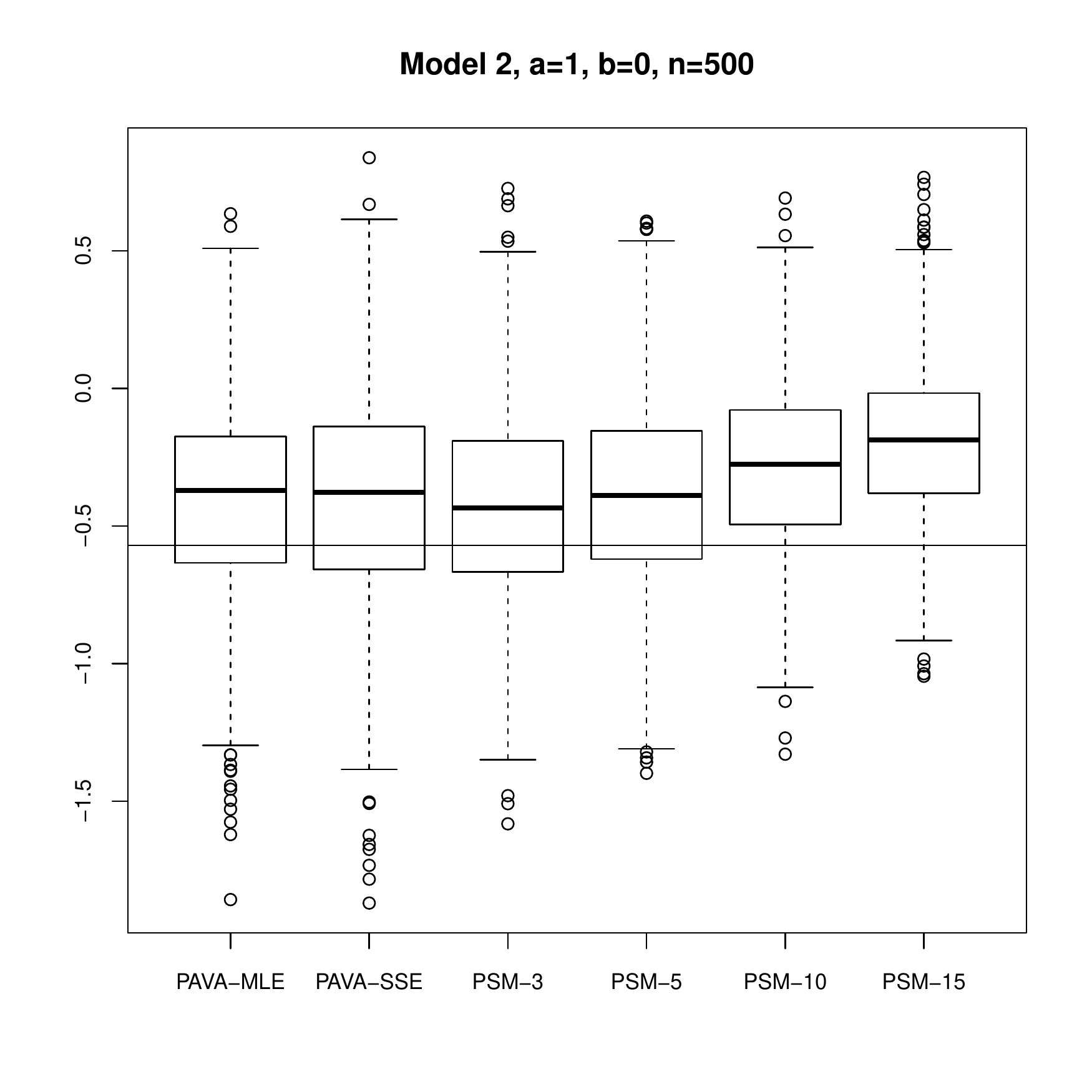} \\
\includegraphics[width=0.45\textwidth, height=0.24\textheight]{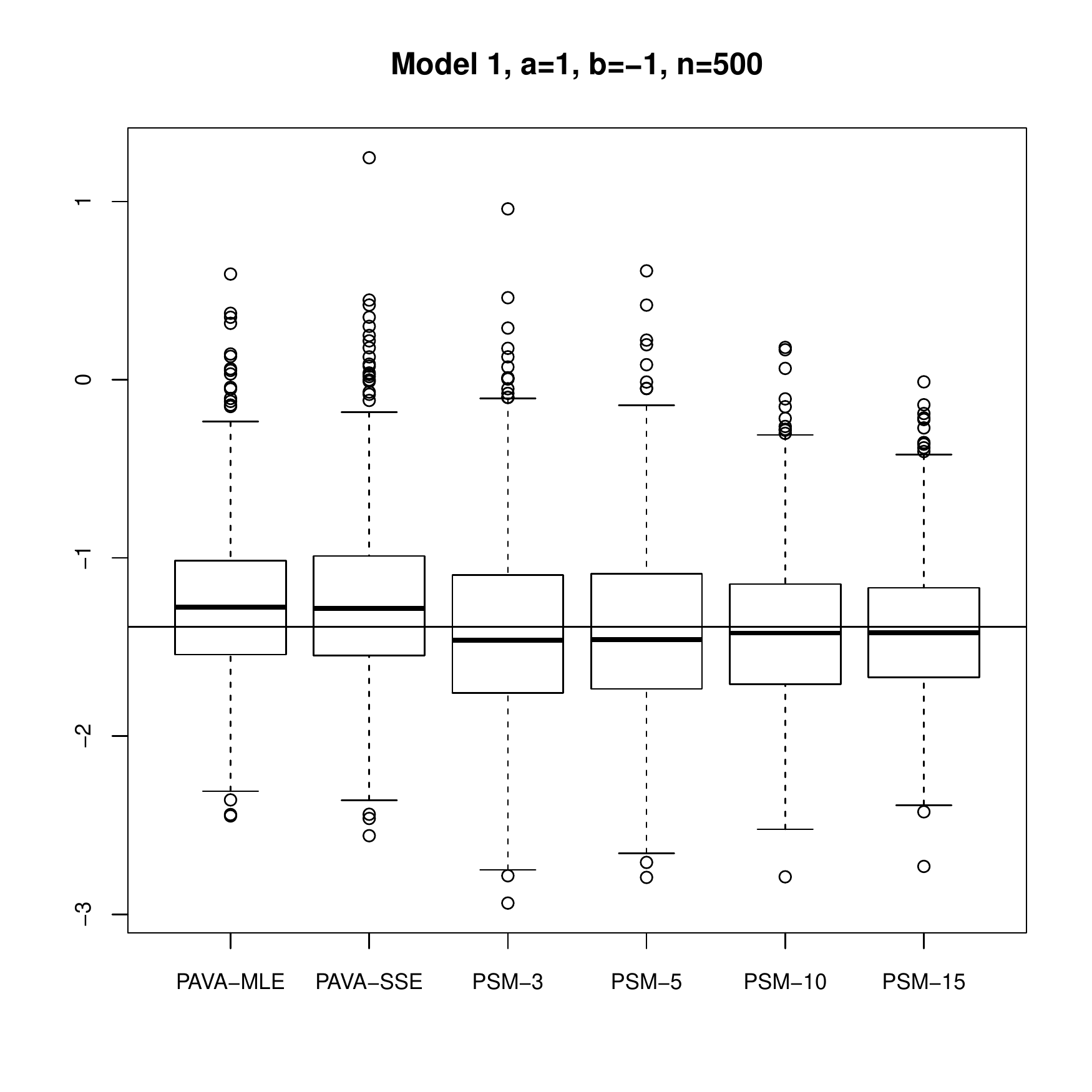}
\includegraphics[width=0.45\textwidth, height=0.24\textheight]{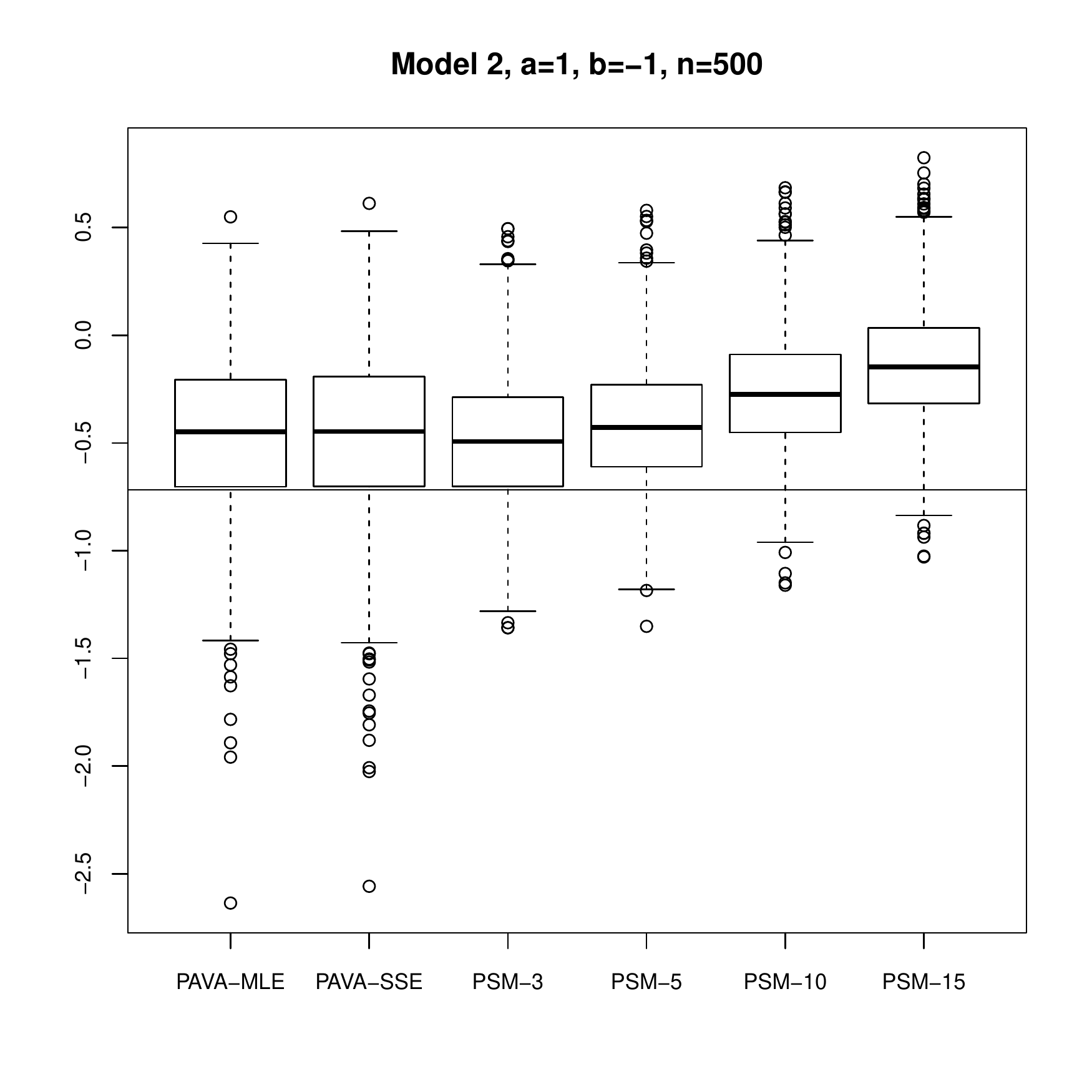} \\
\caption{Boxplots of the seven estimators when $\pi(\cdot)$ is correctly specified,
$a=1$, and $n=500$. As $b$ varies between 1, 0, and -1, the angle between $\beta$
and $\beta$ increases from 0 degrees to 90 degrees.
Model 1: $ h(X_1, X_2) = \cos( X_1+ a X_2)$;
Model 2: $ h(X_1, X_2) = X_1$.
}
\label{qqplot-logistic}
\end{figure}

\begin{figure}[p]
\centering
\includegraphics[width=0.45\textwidth, height=0.24\textheight]{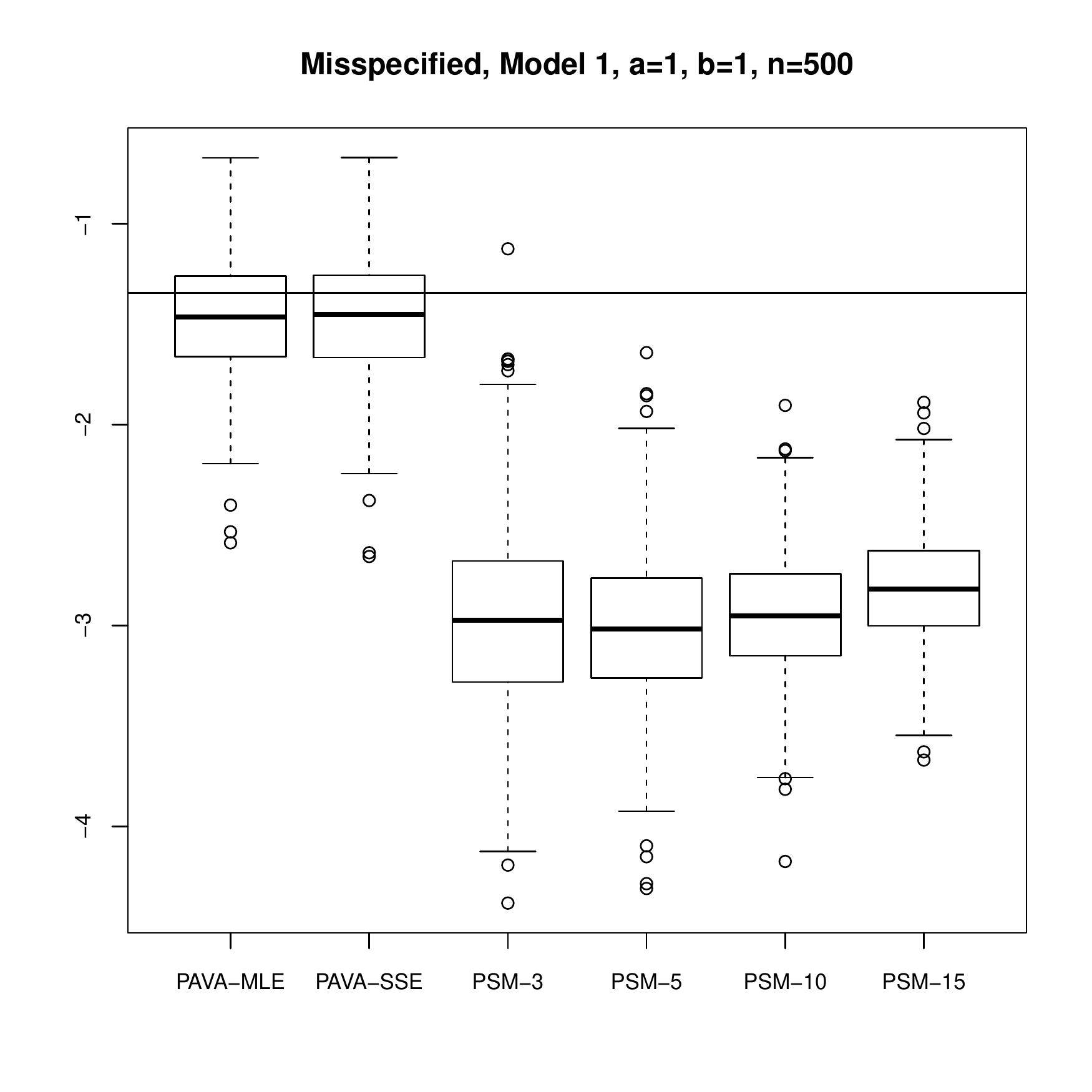}
\includegraphics[width=0.45\textwidth, height=0.24\textheight]{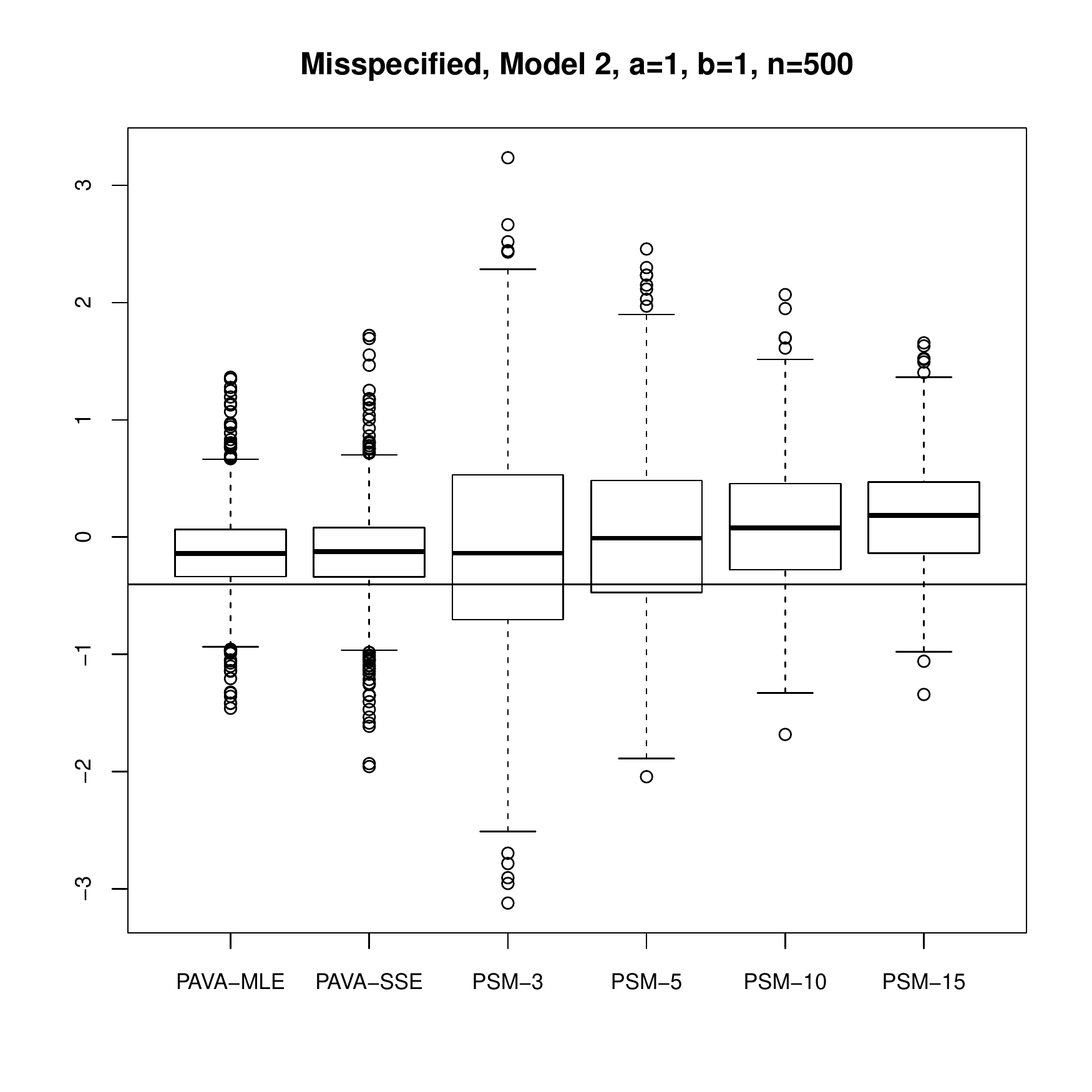} \\
\includegraphics[width=0.45\textwidth, height=0.24\textheight]{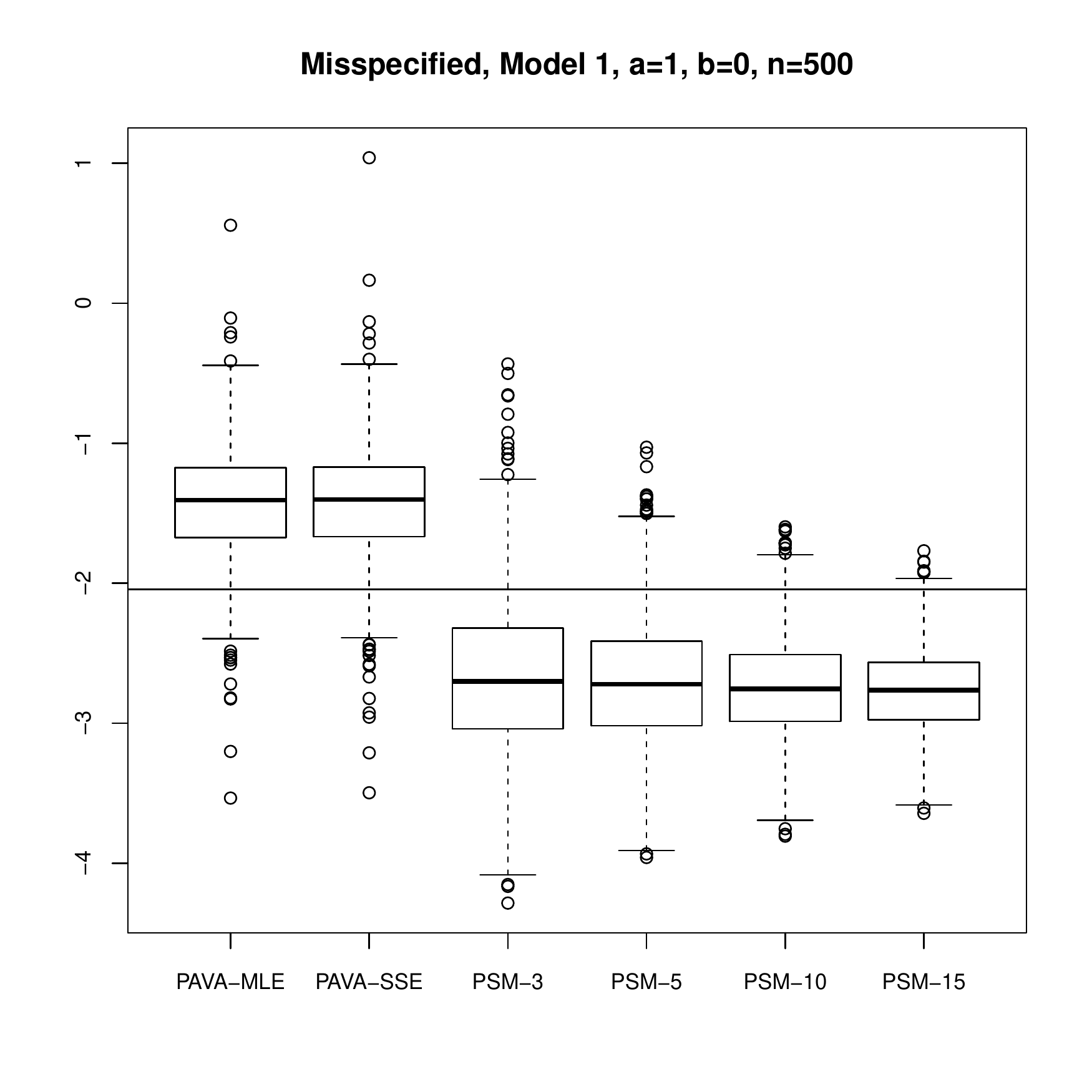}
\includegraphics[width=0.45\textwidth, height=0.24\textheight]{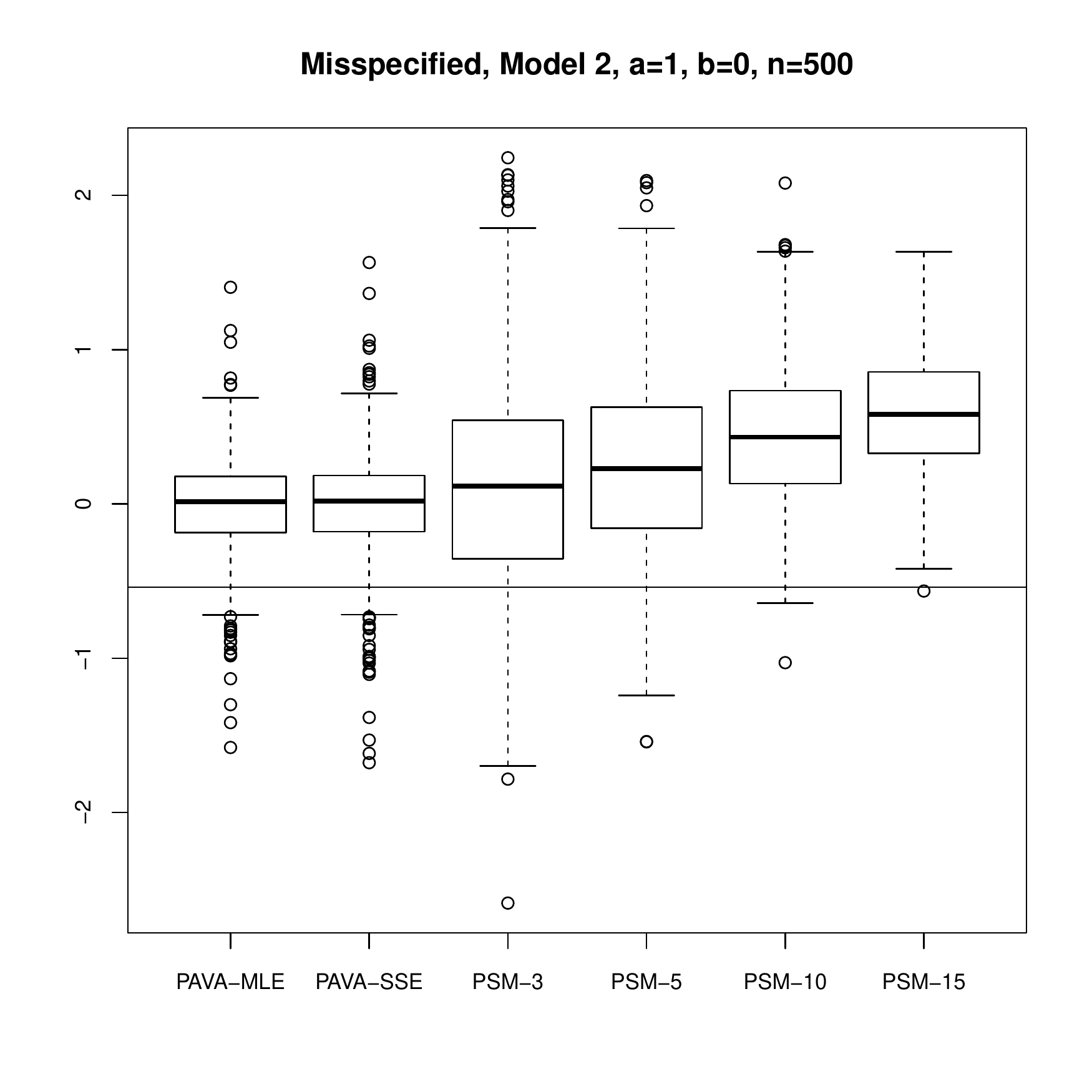} \\
\includegraphics[width=0.45\textwidth, height=0.24\textheight]{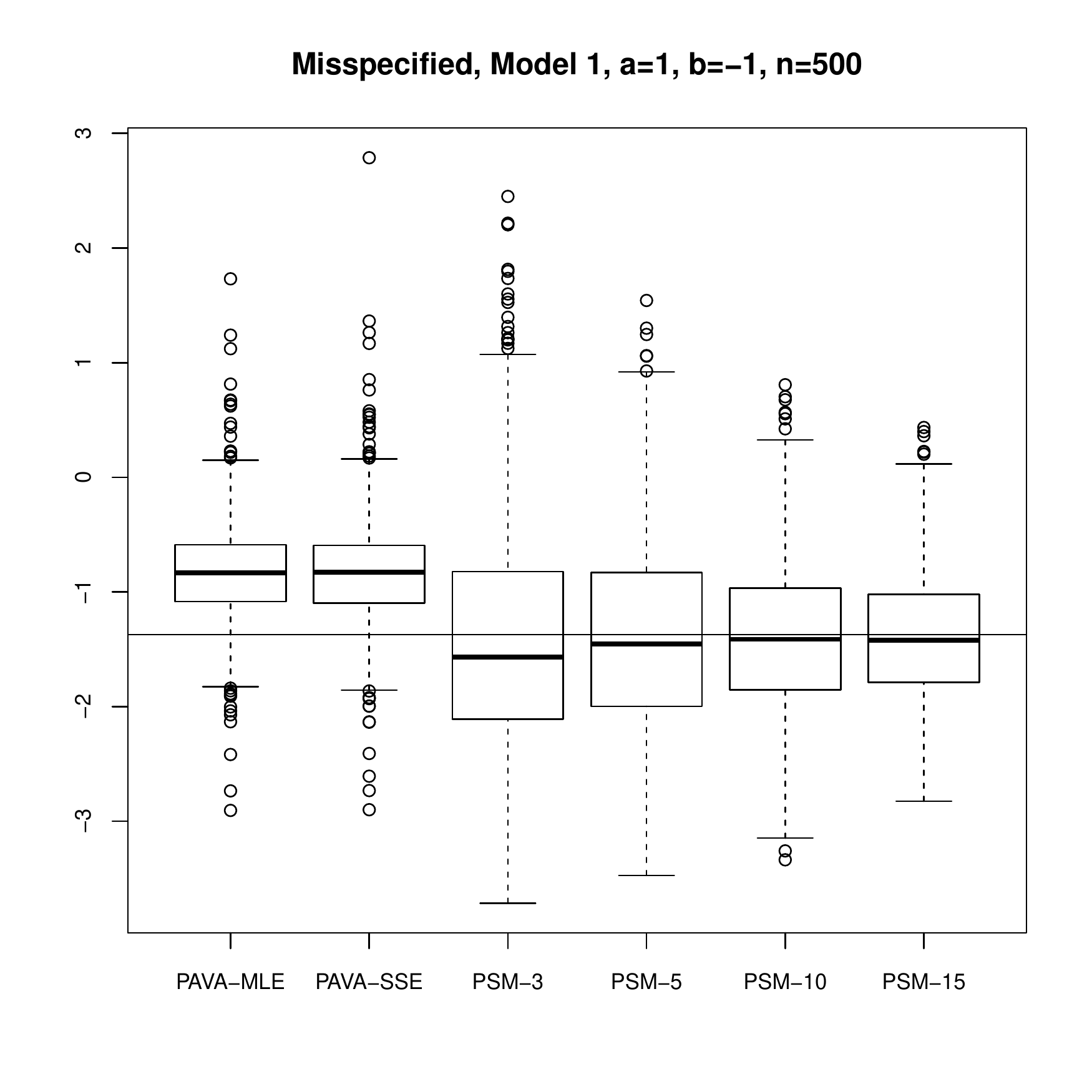}
\includegraphics[width=0.45\textwidth, height=0.24\textheight]{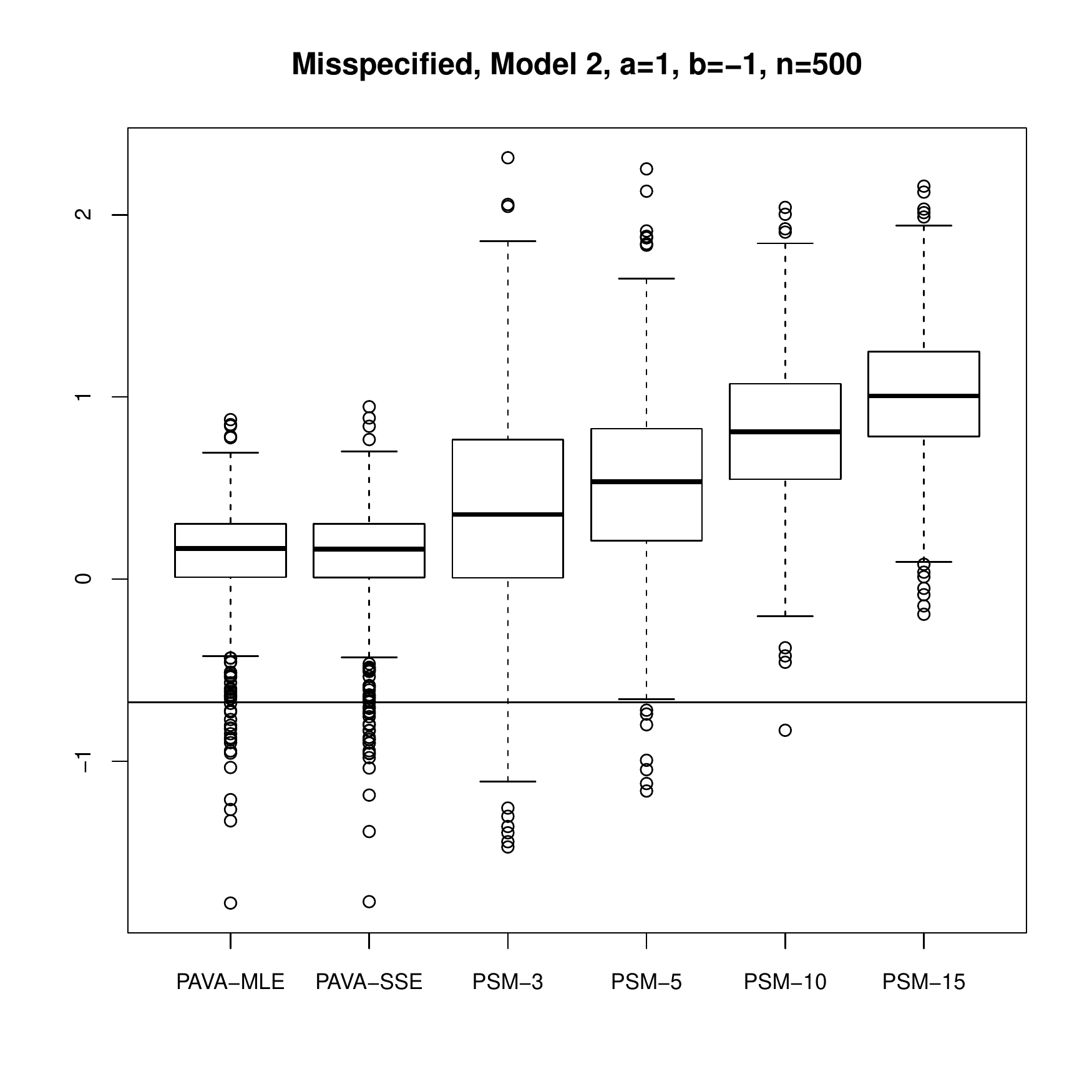} \\
\caption{Boxplots of the seven estimators when $\pi(\cdot)$ is chosen to be the standard normal distribution function,
$a=1$, and $n=500$. As $b$ varies between 1, 0, and -1, the angle between $\beta$
and $\beta$ increases from 0 degrees to 90 degrees.
Model 1: $ h(X_1, X_2) = \cos( X_1+ a X_2)$;
Model 2: $ h(X_1, X_2) = X_1$.
}
\label{qqplot-probit}
\end{figure}

To obtain further insights into the performance of the proposed estimators
and the PSM estimators,
Figures \ref{qqplot-logistic} and \ref{qqplot-probit} display their boxplots in the case of $a=1$;
those of the PARA estimates are excluded as they spread too widely.
As we can see from the boxplots, PAVA-MLE and PAVA-SSE
have smaller biases and better overall performance in most cases
compared with the four PSM estimators.
As the number of matches $M$ increases from 3 to 15,
the PSM estimators exhibit decreasing variances, but
increasing biases except in the case of Model 1, $a=1$, and $b=-1$.
The advantage of the proposed estimators over
the PSM estimators is clarified in Figure \ref{qqplot-probit},
where $\pi(t) $ is chosen to be the standard normal distribution function.

Overall, the proposed PAVA-based estimation method is more reliable
than the PSM method: it usually has smaller RMSEs and biases,
and is not influenced by the tuning parameters.
The performance of the PSM method is strongly influenced by the number of
matches per unit; however, determining the optimal value of $M$ is quite challenging
and no research has been done on this issue.
The direction in which the control response depends on the covariate
does not influence the consistency, but does affect the estimation efficiency of the proposed method:
if it coincides with the direction in which the treatment status depends on the covariate,
the proposed method achieves optimal performance both numerically and theoretically.

\begin{table}[h]
\caption{Simulated biases and RMSEs based on 1000 samples of size $n=500$
when $\pi(\cdot)$ is correctly specified.
Model=1: $h=\cos(X_1+aX_2)$; Model=2: $h=X_1$.
\label{logistic-tab}
}
\centering
\tabcolsep 2pt
\vspace{1ex}
\renewcommand{\arraystretch}{1}
\begin{tabular}{cccc|ccc|cccc}
\hline
Model & $a$ & $b$ & & PAVA-MLE & PAVA-SSE & PARA & PSM-3 & PSM-5 & PSM-10 & PSM-15 \\ \hline
1& 1& 1& Bias & -0.167& -0.158& -0.023& -0.482& -0.646& -0.933& -1.117 \\
1& 1& 1& RMSE & 0.458& 0.468& 1.063& 0.614& 0.744& 0.992& 1.156 \\
1& 1& 0& Bias & 0.125& 0.128& 0.012& -0.200& -0.251& -0.357& -0.424 \\
1& 1& 0& RMSE & 0.407& 0.414& 0.552& 0.463& 0.446& 0.478& 0.507 \\
1& 1& -1& Bias & 0.128& 0.129& -0.038& -0.025& -0.020& -0.017& -0.017 \\
1& 1& -1& RMSE & 0.448& 0.473& 0.584& 0.521& 0.484& 0.429& 0.399 \\\cline{2-11}
1& 2& 1& Bias & -0.599& -0.604& -0.148& -0.901& -1.124& -1.467& -1.639 \\
1& 2& 1& RMSE & 0.936& 0.955& 1.959& 1.039& 1.222& 1.522& 1.677 \\
1& 2& 0& Bias & -0.024& -0.024& -0.010& -0.257& -0.311& -0.399& -0.441 \\
1& 2& 0& RMSE & 0.506& 0.515& 0.878& 0.575& 0.556& 0.571& 0.585 \\
1& 2& -1& Bias & -0.122& -0.141& -0.007& -0.031& -0.043& -0.039& -0.024 \\
1& 2& -1& RMSE & 0.865& 0.895& 1.146& 0.986& 0.906& 0.827& 0.788 \\ \hline
2& 1& 1& Bias & 0.071& 0.080& 0.008& 0.062& 0.086& 0.145& 0.191 \\
2& 1& 1& RMSE & 0.380& 0.444& 0.572& 0.447& 0.403& 0.372& 0.383 \\
2& 1& 0& Bias & 0.159& 0.164& 0.036& 0.138& 0.189& 0.291& 0.381 \\
2& 1& 0& RMSE & 0.380& 0.411& 0.602& 0.391& 0.396& 0.427& 0.478 \\
2& 1& -1& Bias & 0.243& 0.241& 0.017& 0.230& 0.298& 0.452& 0.581 \\
2& 1& -1& RMSE & 0.449& 0.461& 0.812& 0.387& 0.418& 0.531& 0.644 \\\cline{2-11}
2& 2& 1& Bias & -0.319& -0.338& -0.009& -0.308& -0.343& -0.341& -0.299 \\
2& 2& 1& RMSE & 0.553& 0.615& 1.222& 0.561& 0.552& 0.520& 0.479 \\
2& 2& 0& Bias & 0.034& 0.030& 0.002& 0.095& 0.147& 0.262& 0.374 \\
2& 2& 0& RMSE & 0.360& 0.416& 0.566& 0.440& 0.437& 0.459& 0.532 \\
2& 2& -1& Bias & -0.080& -0.103& -0.035& 0.172& 0.239& 0.390& 0.533 \\
2& 2& -1& RMSE & 0.622& 0.662& 0.936& 0.770& 0.739& 0.751& 0.802 \\ \hline
\end{tabular}

\end{table}

\begin{table}[h]
\caption{Simulated biases and RMSEs based on 1000 samples of size $n=500$
when $\pi(\cdot)$ is chosen to be the standard normal distribution function.
Model=1: $h=\cos(X_1+aX_2)$; Model=2: $h=X_1$.
\label{probit-tab}
}
\centering
\tabcolsep 2pt
\vspace{1ex}
\renewcommand{\arraystretch}{1}
\begin{tabular}{cccc|ccc|cccc}
\hline
Model & $a$ & $b$ & & PAVA-MLE & PAVA-SSE & PARA & PSM-3 & PSM-5 & PSM-10 & PSM-15 \\ \hline
1& 1& 1& Bias & -0.114& -0.109& -0.471& -1.624& -1.667& -1.605& -1.474 \\
1& 1& 1& RMSE & 0.309& 0.318& 0.636& 1.686& 1.709& 1.634& 1.501 \\
1& 1& 0& Bias & 0.608& 0.612& 0.278& -0.631& -0.656& -0.693& -0.712 \\
1& 1& 0& RMSE & 0.714& 0.724& 0.722& 0.850& 0.812& 0.784& 0.778 \\
1& 1& -1& Bias & 0.541& 0.535& 0.384& -0.068& -0.019& -0.022& -0.019 \\
1& 1& -1& RMSE & 0.696& 0.708& 0.912& 0.983& 0.847& 0.660& 0.557 \\\cline{2-11}
1& 2& 1& Bias & -0.747& -0.747& -0.930& -1.977& -1.846& -1.414& -0.918 \\
1& 2& 1& RMSE & 0.815& 0.821& 1.069& 2.034& 1.900& 1.495& 1.068 \\
1& 2& 0& Bias & 0.333& 0.333& 0.151& -0.512& -0.500& -0.410& -0.302 \\
1& 2& 0& RMSE & 0.552& 0.560& 1.118& 0.934& 0.837& 0.695& 0.591 \\
1& 2& -1& Bias & -0.521& -0.555& -0.426& -0.177& -0.148& -0.126& -0.135 \\
1& 2& -1& RMSE & 0.971& 0.994& 1.222& 1.880& 1.592& 1.268& 1.115 \\ \hline
2& 1& 1& Bias & 0.266& 0.269& 0.268& 0.312& 0.405& 0.503& 0.582 \\
2& 1& 1& RMSE & 0.459& 0.494& 0.667& 0.955& 0.843& 0.739& 0.741 \\
2& 1& 0& Bias & 0.529& 0.532& 0.554& 0.652& 0.785& 0.974& 1.127 \\
2& 1& 0& RMSE & 0.608& 0.626& 0.829& 0.935& 0.970& 1.069& 1.188 \\
2& 1& -1& Bias & 0.808& 0.804& 0.812& 1.054& 1.207& 1.485& 1.690 \\
2& 1& -1& RMSE & 0.855& 0.853& 1.305& 1.197& 1.300& 1.538& 1.727 \\\cline{2-11}
2& 2& 1& Bias & -0.361& -0.370& -0.199& 0.005& 0.178& 0.660& 1.103 \\
2& 2& 1& RMSE & 0.528& 0.563& 0.640& 0.918& 0.841& 0.939& 1.271 \\
2& 2& 0& Bias & 0.275& 0.291& 0.433& 0.816& 0.986& 1.289& 1.549 \\
2& 2& 0& RMSE & 0.452& 0.520& 0.685& 1.284& 1.307& 1.467& 1.665 \\
2& 2& -1& Bias & -0.280& -0.294& 0.052& 0.941& 1.138& 1.419& 1.629 \\
2& 2& -1& RMSE & 0.650& 0.726& 0.975& 1.708& 1.650& 1.697& 1.816 \\ \hline
\end{tabular}

\end{table}

\section{An application}

For illustration, we apply the proposed PAVA-based estimation method to
data from the National Supported Work (NSW) Demonstration,
which have previously been analyzed by
\cite{LaLonde1986,Dehejia2002}, and \cite{Smith2005}.
The primary parameters of interest in these papers concern
the average treatment effect of a job training program.
We focus on the estimation of the average treatment effect on
the treated, $\tau$. The data consist of 297 treated and 425 untreated observations.
We take earnings in 1978 as the outcome variable of interest ($Y$) and
take age and education as the basic covariates $X_1$ and $X_2$, respectively.
To examine the sensitivity of the proposed method to the model specifications,
we model the propensity score
using a linear logistic model and a quadratic logistic model.

We calculate the point estimates of the seven estimation methods
considered in the previous section.
We also conduct bootstrap sampling with 1000 bootstrap replications from the LaLonde data
to obtain the 2.5\% and 97.5\% sample quantiles,
the mean, and the standard deviation of the resulting point estimates for each method.
Each pair of 2.5\% and 97.5\% sample quantiles constitutes
a 95\% confidence interval of the percentile method.
The analysis results are presented in Table \ref{realdata-tab}.

Under either the linear or quadratic logistic model,
the point estimates of the PAVA-MLE and PAVA-SSE methods
are all around 915 with bootstrap standard deviations of around 500.
This implies that the proposed estimation methods are rather
robust to different model specifications.
As PAVA-SSE makes the weakest model assumption, we believe that
the results of PAVA-SSE should be the most trustable among the seven methods considered here.
The bootstrap means of all methods are around 950
and the point estimates of PSM-15 are about 900 in both cases,
which seemingly provide evidence for the rationality of the PAVA-MLE and PAVA-SSE point estimates.
Although PARA also has bootstrap standard deviations of around 500,
it produces very different point estimates (875.37 and 809.43) in the
linear and quadratic logistic propensity score models.
The PSM method is rather sensitive to the number of matches per unit.
Its point estimate changes from 514.42 to nearly 900 with the linear logistic model,
and varies even more dramatically with the quadratic logistic model.

\begin{figure}
\centering
\includegraphics[width=0.45\textwidth ]{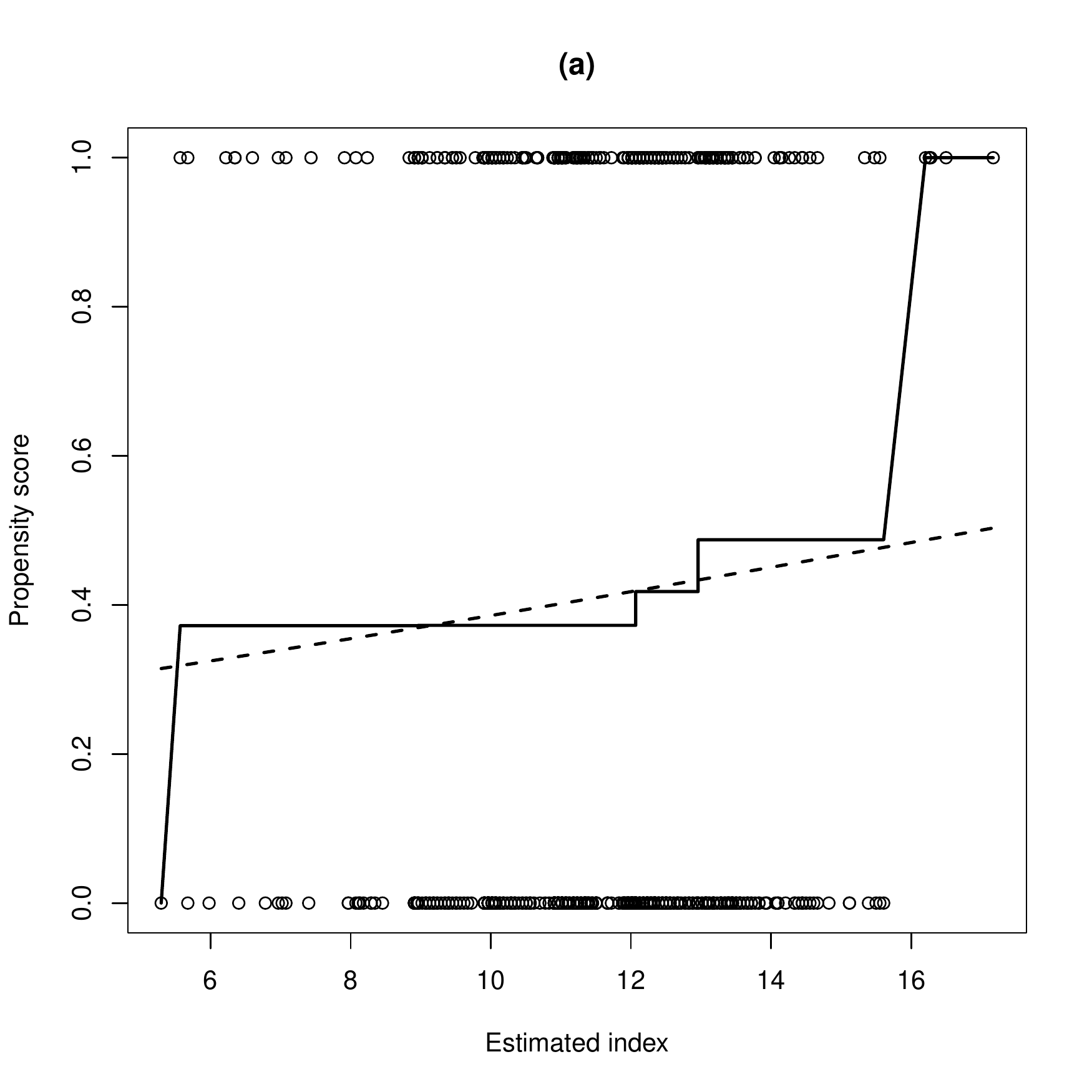}
\includegraphics[width=0.45\textwidth ]{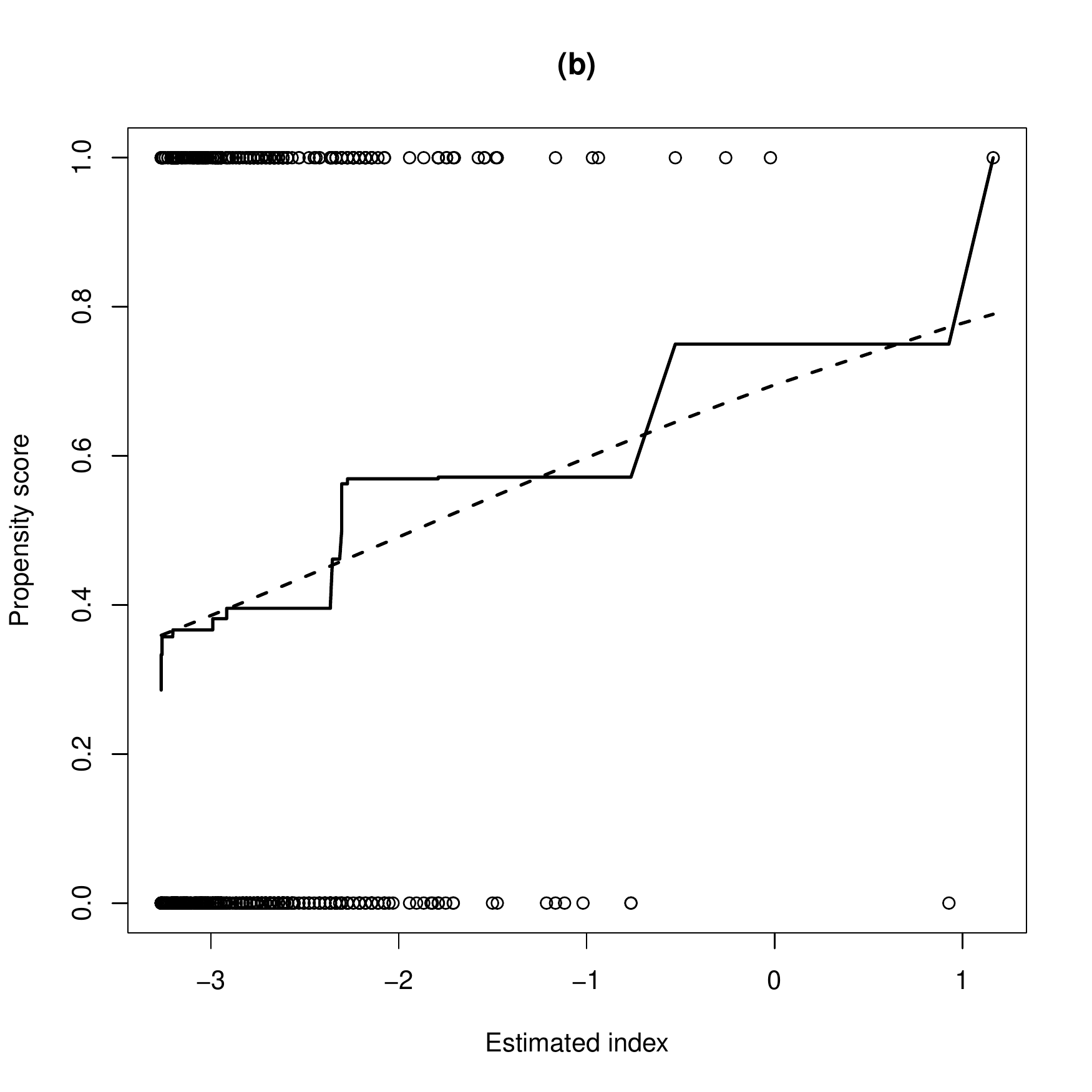} \\
\caption{
Fitted propensity scores versus the estimated index $X^\T \hat \beta$
under a linear (left) and quadratic (right) logistic propensity
score model based on the Lalonde data.
Here, $\hat \beta$ is the MLE of $\beta$ under the corresponding linear logistic model.
Solid line: link function estimated by PAVA;
dashed line: link function set to the logistic function.
}
\label{results-realdata}
\end{figure}

Figure \ref{results-realdata}
displays the fitted propensity scores (versus the estimated index
$X^\T \hat \beta$) using a parametrically logistic model and
the estimations of the semiparametric PAVA method after $\beta$ is replaced by its MLE
under the logistic model.
The parametric propensity score estimates for both the linear and quadratic logistic models apparently form straight lines;
unlike the semiparametric PAVA-based propensity score estimates,
they may not capture local changes in the propensity score.
As the semiparametric method requires fewer model assumptions and is more flexible,
we believe that the semiparametric PAVA-based propensity score estimates
and the corresponding PAVA-MLE estimates are more reliable than
those based on the parametric propensity score estimates,
including PARA and the four PSM methods.
This may explain why the proposed PAVA-based method
is superior to PARA and the four PSM methods.

\begin{table}
\caption{
Estimation results for ATT based on the Lalonde data. 2.5\% and 97.5\% quantiles, Bootstrap mean, and Standard Deviation
denote the corresponding characteristics of the point estimates
based on 1000 bootstrap samples from the Lalonde data.
\label{realdata-tab}
}
\centering
\tabcolsep 3pt
\vspace{1ex}
\renewcommand{\arraystretch}{1}
\begin{tabular}{ cccccccc}
\hline
Methods & PAVA-MLE & PAVA-SSE & PARA & PSM-3 & PSM-5 & PSM-10 & PSM-15 \\ \hline
& \multicolumn{7}{l}{Case (a): $X_1$ and $X_2$ are covariates } \\
Point estimate & 917.33& 911.47& 875.37& 514.42& 714.39& 802.88& 898.99 \\
2.5\% quantile & -41.13& -43.71& -74.50& -156.27& -50.18& -20.66& -12.11 \\
97.5\% quantile & 1857.84& 1914.45& 1815.51& 2202.33& 2086.06& 2009.02& 1966.01 \\
Bootstrap mean & 894.19& 906.39& 861.61& 1071.51& 1033.38& 993.00& 985.47 \\
Standard deviation& 496.56& 495.75& 487.74& 603.06& 558.29& 523.40& 510.66 \\ \hline
& \multicolumn{7}{l}{Case (b): $X_1$, $X_2$, $X_1X_2$, $X_1^2$ and $X_2^2$ are covariates } \\
Point estimate & 913.10& 918.46& 809.43& 88.21& 337.75& 947.25& 903.67 \\
2.5\% quantile & -9.21& -27.03& -119.09& -204.43& -93.13& -40.47& -23.02 \\
97.5\% quantile & 1925.33& 1964.92& 1849.40& 2257.50& 2116.52& 2065.41& 2015.76 \\
Bootstrap mean & 936.21& 950.03& 835.23& 997.70& 997.80& 1019.41& 1028.54 \\
Standard deviation& 507.33& 508.74& 503.50& 645.31& 592.54& 532.30& 515.47 \\ \hline
\end{tabular}

\end{table}

\section{Discussion}

Among many others,
\cite{Abadie2006,Abadie2011,Abadie2012,Abadie2016}
have established different matching methods for causal inference.
Moreover, motivated by the empirical likelihood method
in the presence of auxiliary information and choice-based sampling,
\cite{Hirano2003} proposed an efficient
inverse weighting method using the fully nonparametric estimated propensity score.
Even though their theoretical results are elegant,
the finite-sample performance of their method is unclear.
In practical applications, the dimension of the explanatory variable is high,
and the fully nonparametric estimation of
the regression function may suffer from the curse of dimensionality.
Compared with other efficient estimates, \cite{Hirano2003} stated that
``Which estimators have more attractive finite sample properties, and
which have more attractive computational properties, remain open questions."
The connection between their matching methods and
the probabilistic inverse weighting method, however, is unclear.
In contrast to \cite{Hirano2003}, this paper
has proposed an inverse weighting method that uses
the maximum shape-restricted semiparametric likelihood
estimation of the monotone index propensity score.
Our method is very easy to implement using existing
statistical software in {\tt R}, such as
the {\tt Iso} and {\tt Isotone} packages.
Remarkably, our inverse weighting method is seamlessly
related to the tuning-parameter-free propensity matching method.
Theoretical results show that our estimates can achieve
the semiparametric efficiency lower bound for the average
treatment effect and the average treatment effect for the treated
if the explanatory variable is univariate or the regression function
and propensity score depend
on the explanatory variables in the same direction.
Our results underline the important role played by the propensity score
and the regression function in estimating average causal effects.
In general, the propensity score matching method or the regression function matching method
alone cannot be efficient. An efficient estimation method
should take both of them into consideration \citep{Hu2012}.

\cite{Henmi2004} observed a paradox associated with parameter estimation in the
presence of nuisance parameters. In particular, they found that
the inverse probability weighting estimator with an estimated proposed score
has a smaller asymptotic variance than that derived from the true propensity score.
This paradox was also observed by \cite{Abadie2016}, i.e., matching estimators based on
estimated propensity scores have smaller asymptotic variances than
those based on the true PSM.
Our results further echo this message
that matching estimators based on
the shape-restricted nonparametric MLE
of the propensity score have smaller asymptotic variances
than their counterparts based on the parametric MLE of the
propensity score.

\section*{ Acknowledgements}

%The authors thank the Editor, the Associate Editor,
%and two referees for helpful comments and suggestions
%that have led to significant improvements in the paper.
Yukun Liu's research is supported by
the National Natural Science Foundation of China (12171157),
the State Key Program of National Natural Science Foundation of China (71931004 and 32030063),
and the 111 project (B14019).

\end{document}